\RequirePackage{rotating} 
\documentclass[ijoc,noblindrev]{informs3} 

\OneAndAHalfSpacedXI 

\usepackage{endnotes}
\let\footnote=\endnote
\let\enotesize=\normalsize
\usepackage{natbib}
 \bibpunct[, ]{(}{)}{,}{a}{}{,}%
 \def\bibfont{\footnotesize}%
 \def\BIBand{and}%

\usepackage[utf8]{inputenc}

\usepackage[colorinlistoftodos,prependcaption,textsize=small]{todonotes}
\usepackage{rotating}
\usepackage[shortlabels]{enumitem}
\usepackage{placeins}
\usepackage{mathtools}
\usepackage[noend]{algorithm2e}
\usepackage[english]{babel}
\usepackage[margin=1in]{geometry}
\usepackage{bbding}
\usepackage[shortlabels]{enumitem}
\usepackage[super]{nth}
\usepackage{bm, bbm}
\usepackage{amsmath, amssymb}
\usepackage{mathtools}
\usepackage{subcaption, placeins}
\usepackage{graphicx}
\usepackage[hidelinks]{hyperref}
\usepackage{pgfplots}
\usepackage{ulem}

\usepackage{tikz}
\usetikzlibrary{calc}

\usepackage{subcaption}
\usepackage{booktabs}
\usepackage{multirow} 
\usepackage{placeins}
\usepackage{booktabs}
\usepackage{accents}
\usepackage{graphics}

\usepackage{tikz}
\usepackage{placeins}
\usepackage{float}
\usepackage{graphicx}
\usepackage{accents}
\usepackage{rotating}
\usepackage{graphics}

\usepackage{duckuments}

\usepackage{color}
\definecolor{blue}{rgb}{0,0,1}

\newcommand{\rev}[1]{{\color{black}#1}}

\setlist[itemize]{leftmargin=*, nosep}

\newcommand{\ag}[1]{{\color{black}#1}}

\usepackage{etoolbox}
\newcommand{\zerodisplayskips}{%
  \setlength{\abovedisplayskip}{7pt}%
  \setlength{\belowdisplayskip}{5pt}%
  \setlength{\abovedisplayshortskip}{0pt}%
  \setlength{\belowdisplayshortskip}{0pt}}
\appto{\normalsize}{\zerodisplayskips}
\appto{\small}{\zerodisplayskips}
\appto{\footnotesize}{\zerodisplayskips}
\usepackage{xfrac}

\TheoremsNumberedThrough     
\ECRepeatTheorems

\EquationsNumberedThrough    

\begin{document}


\RUNAUTHOR{Ryan Cory-Wright and Andr\'es G\'omez} 

\RUNTITLE{Efficient Cross-Validation for Sparse Linear Regression}

\TITLE{Efficient Cross-Validation \\ for Sparse Linear Regression}

\ARTICLEAUTHORS{%
\AUTHOR{Ryan Cory-Wright}
\AFF{Department of Analytics, Marketing and Operations, Imperial Business School, London, UK\\
ORCID: \href{https://orcid.org/0000-0002-4485-0619}{$0000$-$0002$-$4485$-$0619$}\\ \EMAIL{r.cory-wright@imperial.ac.uk}} 
\AUTHOR{Andr{\'e}s G{\'o}mez}
\AFF{Department of Industrial and Systems Engineering, Viterbi School of Engineering, University of Southern California, CA\\
ORCID: \href{https://orcid.org/0000-0003-3668-0653}{$0000$-$0003$-$3668$-$0653$}\\
\EMAIL{gomezand@usc.edu}}
} 

\ABSTRACT{%
Given a high-dimensional covariate matrix and a response vector, ridge-regularized sparse linear regression selects a subset of features that explains the relationship between covariates and the response in an interpretable manner. To 
{\color{black}choose hyperparameters that control the sparsity level and amount of regularization, practitioners commonly use} $k$-fold cross-validation. However, cross-validation substantially increases the computational cost of sparse regression as it requires solving many mixed-integer optimization problems (MIOs) for each hyperparameter combination. To {\color{black}address this computational burden}, we {\color{black}derive} computationally tractable relaxations of {\color{black}the} $k$-fold cross-validation {\color{black}loss}, facilitating hyperparameter selection {\color{black}while} solving $50$--$80\%$ fewer MIOs in practice. {\color{black}Our computational results demonstrate, across eleven real-world UCI datasets, that MIO-based cross-validation can be competitive with mature software packages such as glmnet and L0Learn.}
}%


\maketitle
\vspace{-10mm}

\section{Introduction}\label{sec:intro}
Over the past fifteen years, Moore's law has spurred an 
explosion of high-dimensional datasets for scientific discovery across multiple fields \citep[][]{mcafee2012big}. These datasets often consist of a design matrix $\bm{X} \in \mathbb{R}^{n \times p}$ of explanatory variables and an output vector $\bm{y} \in \mathbb{R}^n$ of response variables. Accordingly, practitioners often aim to explain the response variables linearly via the equation $\bm{y}=\bm{X}\bm{\beta}+\bm{\epsilon}$. {\color{black}Using this equation, the vector of regression coefficients $\bm{\beta} \in \mathbb{R}^p$ is inferred by minimizing the least squares (LS) error of the residuals $\bm{\epsilon}$}. 

Despite its computational efficiency, LS regression exhibits two practical limitations. First, when $p \gg n$, there is not enough data to accurately infer $\bm{\beta}$ via LS, and LS regression generates estimators which perform poorly out-of-sample due to {\color{black}the} curse of dimensionality \citep[][]{buhlmann2011statistics, gamarnik2022sparse}.
Second, LS regression generically selects every feature, including irrelevant ones. This is a significant challenge when regression coefficients are used for high-stakes decision-making tasks and non-zero coefficients {\color{black}guide decisions}.

To tackle the challenges of dimensionality and false discovery, sparse learning has emerged as a popular methodology for explaining the relationship between inputs {\color{black}$\bm{X}$} and outputs {\color{black}$\bm{y}$}. 
A popular sparse learning model is ridge-regularized sparse regression, which admits the formulation \citep[][]{bertsimas2020sparse,xie2020scalable, hastie2020best,  atamturk2020safe, kenney2021mip, hazimeh2021sparse, liu2023okridge}
\begin{align}\label{eqn:l0l2_training}
    \min_{\bm{\beta} \in \mathbb{R}^p} \quad &  \frac{\gamma}{2} \Vert \bm{\beta}\Vert_2^2+\Vert \bm{y}-\bm{X}\bm{\beta}\Vert_2^2  \quad \text{s.t.} \quad \Vert \bm{\beta}\Vert_0 \leq \tau,
\end{align} 
where $\tau \in \{1, \ldots, p\}$ and $\gamma >0$ are hyperparameters that respectively {\color{black}control the sparsity of $\bm{\beta}$ and the amount of $\ell_2^2$ regularization \citep[cf.][]{xu2008robust, bertsimas2018characterization}, and we assume 
that $\bm{X}$ {\color{black}and} $\bm{y}$ have undergone standard preprocessing so that $\bm{y}$ is a zero-mean vector and $\bm{X}$ has zero-mean, unit-variance columns, meaning $\gamma$ penalizes each feature equally.

{\color{black}Problem \eqref{eqn:l0l2_training} belongs to the class of sparse quadratic optimization problems, which are known to be NP-hard in general \citep{natarajan1995sparse, chen2019approximation} and computationally challenging in practice} and early mixed-integer formulations could not scale to problems with thousands of features \citep[][]{hastie2020best}. In a more positive direction, by developing and exploiting tight conic relaxations of appropriate substructures of \eqref{eqn:l0l2_training}, e.g., the perspective relaxation \citep{ceria1999convex, stubbs1999branch, gunluk2010perspective}
, more {\color{black}modern} mixed-integer {\color{black}approaches which combine branch-and-bound with cutting planes, presolve, heuristics, and other algorithmic enhancements} \citep{hazimeh2021sparse}{\color{black}, can} scale to larger instances with thousands of 
features. {\color{black}We refer to \cite{bertsimas2021unified, atamturk2019rank} for reviews of perspective and related relaxations.}

To be sure, the aforementioned works solve \eqref{eqn:l0l2_training} rapidly. Unfortunately, they do not address arguably the most significant difficulty {\color{black}in performing sparse regression}. The hyperparameters $(\tau, \gamma)$ 
{\color{black}are not known to the decision-maker ahead of time, as is often assumed in the literature for convenience.}
Rather, they must be selected by {\color{black}the decision-maker}, which is potentially much more challenging than solving \eqref{eqn:l0l2_training} for a single value of $(\tau, \gamma)$ \citep[]{hansen1992new}. Indeed, selecting $(\tau, \gamma)$ typically involves minimizing a validation metric over a grid of values, which is computationally expensive \citep[]{larochelle2007empirical}. 

{\color{black} Perhaps the most popular validation {\color{black}procedure} is hold-out \citep{hastie2009elements}, where one omits a portion of the data when training the model {\color{black}(i.e., solving an instance of \eqref{eqn:l0l2_training})} and then evaluates performance on this hold-out set as a proxy for the model's test set performance. However, hold-out validation is sometimes called a high-variance approach \citep{hastie2009elements}, because the validation errors can vary significantly depending on the hold-out set selected. 

To reduce the variance in this procedure, a number of authors have proposed {\color{black}what we call \textit{the cross-validation paradigm}. }
{\color{black}Early iterations of this paradigm, as reviewed by} \citet{stone1978cross}, {\color{black}suggest solving Problem \eqref{eqn:l0l2_training} a total of $n$ times, each time leaving out a single data point $i \in \rev{\{1,\dots,n\}}$, and estimating out-of-sample performance via the average prediction error of each estimator on its left-out observation.} This approach is known as leave-one-out cross-validation (LOOCV). 

A popular variant of LOOCV, known as $k$-fold cross-validation, {\color{black}involves} removing subsets of $n/k$ data points at a time and breaking the data into $k$ folds in total, which significantly reduces the computational burden of cross-validation while having less variance than a hold-out approach \citep{burman1989comparative, arlot2010survey}. However, even $k$-fold cross-validation may be prohibitive in the case of MIOs such as \eqref{eqn:l0l2_training}. Indeed, as identified by \cite{hastie2020best}, with a time limit of {\color{black}30} minutes per MIO, using 10-fold cross-validation to choose between subset sizes $\tau=1,\dots,50$ in an instance of {\color{black}sparse linear regression} with $p=100$ and $n=500$ requires {\color{black}a cumulative time budget of } {\color{black}250} hours {\color{black}across all of the MIOs}.} 

For sparse regression, {\color{black}given a partition $\mathcal{N}_1,\dots, \mathcal{N}_k$ of {\color{black}$[n]:=\{1, \ldots, n\}$,}  performing {\color{black}$k$-fold cross-validation} corresponds to selecting hyperparameters $\gamma, \tau$ which minimize the function
{\color{black}
\begin{equation}\label{prob:upperlevelofv_validation}
	h(\gamma,\tau) = \frac{1}{\color{black}n}\sum_{j=1}^k\sum_{i\in \mathcal{N}_j} (y_i -\bm{x_i}^\top \bm{\beta}^{(\mathcal{N}_j)}(\gamma,\tau))^2
\end{equation}{\color{black}
where $\bm{\beta}^{(\mathcal{N}_j)}(\gamma,\tau)$ denotes an optimal solution to the following lower-level problem for any $\mathcal{N}_j$:
\begin{equation}
	\bm{\beta}^{(\mathcal{N}_j)}(\gamma,\tau) \in \argmin_{\bm{\beta} \in \mathbb{R}^p} \ \frac{\gamma}{2}\Vert \bm{\beta}\Vert_2^2 +\Vert \bm{y}^{(\mathcal{N}_j)}-\bm{X}^{(\mathcal{N}_j)}\bm{\beta}\Vert_2^2\quad \text{s.t.}\quad \Vert \bm{\beta}\Vert_0 \leq \tau, \label{eqn:lowerlevel1}
\end{equation}}
$\gamma>0$ is a hyperparameter, $\tau$ is a sparsity budget, $\bm{X}^{(\mathcal{N}_j)}, \bm{y}^{(\mathcal{N}_j)}$ denote the dataset with the data in $\mathcal{N}_j$ removed, and we take $\bm{\beta}^{(\mathcal{N}_j)}(\gamma, \tau)$ to be unique for a given $\tau, \gamma$ for convenience}\footnote{This assumption seems plausible, as the training objective is strongly convex for a fixed binary support vector, and therefore for each binary support vector there is indeed a unique solution. }. 
In words, $h(\gamma, \tau)$ denotes the average prediction error on each left-out fold for a sparse regressor with hyperparameters $(\gamma, \tau)$ trained on the remaining folds.

{\color{black}We remark that if all sets $\mathcal{N}_j$ are taken to be singletons and $k=n$, minimizing $h$ corresponds to LOOCV. Moreover, if $k=2$ and the term with $j=2$ is removed from $h$, optimizing $h$ reduces to minimizing the hold-out error.}
After selecting $(\gamma, \tau)$, practitioners usually train a final model on the entire dataset, by solving Problem \eqref{eqn:l0l2_training} with the selected hyperparameter combination. 

\paragraph{Our Approach:} We propose techniques for obtaining strong bounds on validation metrics in polynomial time and leverage these bounds to design algorithms for minimizing the {\color{black}cross-validation} error in Sections \ref{sec:lowerbounds} {\color{black}and} \ref{sec:corddescent}. By performing a perturbation analysis of perspective relaxations of sparse regression problems, we construct convex relaxations of the $k$-fold cross-validation error, which allows us to minimize it without explicitly solving MIOs at each data fold and for each hyperparameter combination. This results in a branch-and-bound algorithm for hyperparameter selection that is substantially more efficient than state-of-the-art {\color{black}methods such as} grid search. \color{black}As an aside, we remark that as cross-validation is more general than hold-out validation, our convex relaxations can be generalized immediately to the hold-out case.

{\color{black}In numerical experiments (Section \ref{sec:numres}), we assess the impact of our contributions. We observe on real UCI datasets that our branch-and-bound scheme reduces the number of MIOs that need to be solved by an average of $50\%$–{\color{black}$80\%$}. Further, we leverage our branch-and-bound scheme to design a cyclic \rev{alternating minimization} scheme that iteratively minimizes $\tau$ and $\gamma$. {\color{black}In experiments with UCI machine learning datasets, it performs {\color{black}competitively} to the glmnet \citep{friedman2010regularization}, MCP \citep{zhang2010nearly} and L0Learn \citep{hazimeh2022l0learn} software packages in terms of solution quality, {\color{black}resulting in the best overall performance in 3/11 datasets tested and within 5\% of the best package in another 5/11 datasets.}

\subsection{Literature Review}
Our work falls at the intersection of three areas of the optimization literature: 
{\color{black}(i)} hyperparameter selection techniques for optimizing the performance of a machine learning model by selecting hyperparameters that perform well on a validation set, {\color{black}(ii)} bilevel approaches that reformulate and solve hyperparameter selection problems as bilevel problems, and {\color{black}(iii)} perspective reformulation techniques for mixed-integer problems with logical constraints{\color{black}, as discussed above.} To put our contributions into context, we now review {\color{black}the two remaining} areas of the literature.

\paragraph{Hyperparameter Selection Techniques for Machine Learning Problems:} 
A wide variety of hyperparameter selection techniques have been proposed for machine learning problems such as sparse regression, 
{\color{black}including grid search \citep[][]{larochelle2007empirical} as reviewed in Section \ref{sec:intro},}
{\color{black}and random search \citep[cf.][]{bergstra2012random}.} 
In random search, we let $\mathcal{L}$ be a random sample from a space of valid hyperparameters, e.g., a uniform distribution over $[10^{-4}, 10^4] \times [p]$ for sparse regression. Remarkably, in settings with many hyperparameters, random search usually outperforms grid search for a given budget on the number of training problems that can be solved, because validation functions often have a lower effective dimension than the number of hyperparameters present in the model \citep[][]{bergstra2012random}. 
However, grid search remains competitive for problems with a small number of hyperparameters, such as sparse regression. 

The modern era of hyperparameter selection strategies was ushered in by the increasing prominence of deep learning methods in applications from voice recognition to drug discovery \citep[see][for a review]{lecun2015deep}. 
The volume of data available and {\color{black}the} number of hyperparameters {\color{black}to be selected have} challenged the aforementioned methods and led to new techniques, including evolutionary strategies
, Bayesian optimization techniques \citep{ frazier2018tutorial} and bandit methods \rev{\citep{falkner2018bohb}}. 
However, in sparse regression problems where we aim to optimize two hyperparameters, these methods are {\color{black}effectively equivalent} to grid or random search. Further, none of these approaches provide locally optimal hyperparameter combinations with respect to {\color{black}the LOOCV error}, which suggests there is room for improvement upon the state-of-the-art in sparse regression. 

We point out that current approaches for hyperparameter selection are similar to existing methods for multi-objective mixed-integer optimization. While there has been recent progress in improving multi-objective algorithms for mixed-integer linear programs \citep{lokman2013finding,stidsen2014branch}, a direct application of these methods might be unnecessarily expensive. Indeed, these approaches seek to compute the efficient frontier \citep{boland2015criterion,boland2015criterion2}{\color{black}, i.e.,} solving problems for all possible values of the regularization parameter. In contrast, we are interested in only the combination of parameters that optimize a well-defined metric (e.g., \ag{\color{black}the cross-validation error}). 
}

\paragraph{Bilevel Optimization for Hyperparameter Selection:} In a complementary direction, several authors have proposed selecting hyperparameters via bilevel optimization \citep[see][for a general theory]{beck2021gentle}, since \cite{bennett2006model} 
recognized that cross-validation is a special case of bilevel optimization. 
Therefore, {\color{black}in principle, we could} minimize the {\color{black}cross-validation} error in sparse regression by invoking bilevel techniques. Unfortunately, this approach seems intractable in both theory and practice \citep[][]{ben1990computational, hansen1992new}. Indeed, standard bilevel approaches such as dualizing the lower-level problem are challenging to apply in our context because our lower-level problems are non-convex and cannot easily be dualized.

Although {\color{black}bilevel hyperparameter optimization is slow in its original implementation}, several authors have proposed making {\color{black}it} more tractable by combining {\color{black}it} with {\color{black}efficient} modeling paradigms to obtain locally optimal sets of hyperparameters. Among others, \cite{sinha2020gradient} recommend taking a gradient-based approximation of the lower-level problem and thereby reducing the bilevel problem to a single-level problem, \cite{okuno2021lp} advocate selecting hyperparameters by solving the KKT conditions of a bilevel problem, and \cite{ye2022difference} propose solving bilevel hyperparameter problems via difference-of-convex methods to obtain a stationary point.

Specializing our review to regression, three works aim to optimize the performance of regression models on a validation metric. First, \cite{takano2020best} propose optimizing the $k$-fold validation loss, assuming all folds share the same support. Unfortunately, although their assumption improves their method's tractability, it may lead to subpar statistical performance {\color{black}because using the same set of non-zero regressors for all folds shares information between the folds}. Second, \citet{stephenson2021can} propose first-order methods for minimizing the leave-one-out error in ridge regression problems (without sparsity constraints). However, it is unclear how to generalize their approach to settings with sparsity constraints. Finally, perhaps closest to our work, \cite{kenney2021mip} propose a bisection algorithm for selecting the optimal sparsity parameter in a sparse regression problem by approximately minimizing the $k$-fold cross-validation error. It is, however, worth noting that this approach is not guaranteed to converge to an optimal sparsity parameter with respect to the $k$-fold error, because it does not develop lower bounds on the $k$-fold error. 

\subsection{Structure}

The rest of the paper is laid out as follows:
\begin{itemize}\color{black}
\item In Section \ref{sec:lowerbounds}, we observe that validation metrics are potentially expensive to evaluate, because they involve solving up to {\color{black}$k$} MIOs (in the $k$-fold case), and accordingly develop tractable lower and upper bounds that can be computed without solving any MIOs.
\item In Section \ref{sec:corddescent}, we propose an efficient \rev{alternating minimization} scheme for identifying locally optimal hyperparameters with respect to the validation error. Specifically, in Section \ref{ssec:parametric}, we develop an efficient scheme for minimizing the {\color{black}cross-validation} error with respect to $\tau$, and in Section \ref{ssec:parametric2}, we propose a scheme for optimizing with respect to $\gamma$.
\item Finally, in Section \ref{sec:numres}, we benchmark our proposed approaches on real {\color{black}UCI} datasets, testing the proposed approach against state-of-the-art software packages for sparse regression. 

\end{itemize}

\subsection*{Notation} We let non-boldface characters such as $b$ denote scalars, lowercase bold-faced characters such as $\bm{x}$ denote vectors, uppercase bold-faced characters such as $\bm{A}$ denote matrices, and calligraphic uppercase characters such as $\mathcal{Z}$ denote sets. We let $[n]$ denote the running set of indices $\{1, \dots, n\}$, and $\Vert \bm{x}\Vert_0:=|\{j: x_j \neq 0\}|$ denote the $\ell_0$ pseudo-norm, i.e., the number of non-zero entries in $\bm{x}$. Finally, we let $\bm{e}$ denote the vector of ones, and $\bm{0}$ denote the vector of all zeros. 


Furthermore, we consistently use notation commonly found in the supervised learning literature. We consider a setting where we observe covariates $\bm{X}:=(\bm{x}_1{\color{black}^\top}; \ldots; \bm{x}_n{\color{black}^\top}) \in \mathbb{R}^{n \times p}$ and response data $\bm{y}:=(y_1, \ldots y_n) \in \mathbb{R}^n$. {\color{black}With this notation, the $i$th row of $\bm{X}$ is denoted by $\bm{x}_i^\top \in \mathbb{R}^p$. }
We say that $(\bm{X}, \bm{y})$ is a training set, and let $\bm{\beta}$ denote a regressor fitted on this training set. In cross-validation, we are also interested in the behavior of $\bm{\beta}$ after {\color{black}leaving out portions of} the training set. We let $(\bm{X}^{(i)}, \bm{y}^{(i)})$ denote the training set with the $i$th data point left out, and denote by $\bm{\beta}^{(i)}$ the regressor obtained after leaving out the $i$th point. {\color{black}Similarly, given a partition $\mathcal{N}_1,\dots,\mathcal{N}_k$  of $[n]$ and $j\in [k]$, we let $(\bm{X}^{(\mathcal{N}_j)}, \bm{y}^{(\mathcal{N}_j)})$ denote the training set with the $j$th fold left out, and $\bm{\beta}^{(\mathcal{N}_j)}$ be the associated regressor.} 

\section{Convex Relaxations of $k$-fold Cross-Validation Error}\label{sec:lowerbounds}

In this section, we develop tractable upper and lower approximations of the $k$-fold cross-validation error of a sparse regression model, which can be evaluated at a given $(\gamma, \tau)$ without solving any MIOs. From a theoretical perspective, one of our main contributions is that, given $\bm{x}\in \mathbb{R}^p$, we show how to construct bounds $\underline \xi, \overline\xi$ such that $\underline \xi\leq \bm{x^\top\beta}^{(\mathcal{N}_j)}\leq \overline \xi$, which we can use to infer out-of-sample predictions. In particular, we leverage this insight to bound from above and below the function:
\begin{align}
	&h(\gamma, \tau)=1/n \sum_{j=1}^k h_j(\gamma, \tau)=1/n\sum_{j=1}^k\sum_{i \in \mathcal{N}_j}\left(y_i -\bm{x_i}^\top \bm{\beta}^{(\mathcal{N}_j)}(\gamma,\tau)\right)^2,\label{eqn:crossValError}
\end{align} 
{\color{black}i.e.,} the $k$-fold cross-validation error. {\color{black}Note that \eqref{eqn:crossValError} is a restatement of \eqref{prob:upperlevelofv_validation}{\color{black}, and $h_j(\gamma, \tau):=\sum_{i \in \mathcal{N}_j}\left(y_i -\bm{x_i}^\top \bm{\beta}^{(\mathcal{N}_j)}(\gamma, \tau)\right)^2$.}
\subsection{Bounds on the Prediction Spread}\label{sec:predictionBound}

Given any \rev{$0<\epsilon\leq \gamma$}, it is well-known that Problem~\eqref{eqn:l0l2_training} admits the conic quadratic relaxation:
\begin{align}\label{eq:persp}
\zeta_{\text{persp}}=\min_{\bm{\beta} \in \mathbb{R}^p,\bm{z}\in [0,1]^p} \quad  \Vert \bm{y}-\bm{X}\bm{\beta}\Vert_2^2 +\rev{\frac{\gamma-\epsilon}{2}} \sum_{i=1}^p \frac{\beta_i^2}{z_i}+\rev{\frac{\epsilon}{2} \sum_{i=1}^p \beta_i^2} \quad \text{s.t.} \quad &\sum_{i=1}^pz_i \leq \tau,
\end{align}
which is also known as the {perspective relaxation} \citep[][]{ceria1999convex, xie2020scalable}. \rev{In the formulation, $\epsilon$ is presumed to be a small number, and the associated regularization term is added to ensure strong convexity. Note that here and throughout the rest of the paper, there is an implicit dependence of $\zeta_{\text{persp}}$ on $\epsilon$.} If integrality constraints $\bm{z}\in \{0,1\}^p$ are added to \eqref{eq:persp}, then the resulting mixed-integer optimization problem (MIO) is a reformulation of \eqref{eqn:l0l2_training}, where the logical constraints {\color{black}$\beta_i=0$ if $z_i= 0 \ \forall i \in [p]$} are implicitly imposed via the domain of the perspective function $\beta_i^2/z_i$. Moreover, the optimal objective $\zeta_{\text{persp}}$ of \eqref{eq:persp} often provides tight lower bounds on the objective value of \eqref{eqn:l0l2_training} \citep[]{pilanci2015sparse,bertsimas2020sparse, askari2022approximation}, and the optimal solution $\bm{\beta}_{\text{persp}}^\star$ is often a good estimator in its own right. As we establish in our main results, the perspective relaxation can also be used to obtain accurate approximations of and bounds on the $k$-fold cross-validation error.

\ag{Our first main result (Theorem~\ref{prop:kfoldspread}) reveals that any optimal solution of \eqref{eqn:l0l2_training} lies in an ellipsoid centered at its continuous (perspective) relaxation, and whose radius depends on the duality gap: 
}

\begin{theorem}\label{prop:kfoldspread}
	Given any bound
	\begin{equation}\label{eq:MIPII}
		\bar u\geq \min_{\bm{\beta}\in \mathbb{R}^p}\; \|\bm{X}\bm{\beta}-\bm{y}\|_2^2+\frac{\gamma}{2}\|\bm{\beta}\|_2^2
		\text{\rm\ s.t. }\; \|\bm{\beta}\|_0\leq \tau,
	\end{equation} 
	{\color{black}and any  $0 \leq \epsilon < \gamma$} the inequality
	\begin{align}\label{ineq:perspbound}
		(\bm{\beta^\star_{\text{persp}}}- \bm{\beta}^\star_{\text{MIO}})^\top \left(\bm{X}^\top \bm{X}\rev{+\frac{\epsilon}{2}\mathbb{I}}\right) (\bm{\beta^\star_{\text{persp}}}- \bm{\beta}^\star_{\text{MIO}}) \leq (\bar{u}-\zeta_{\text{persp}})
	\end{align}
	holds, where $\bm{\beta}_{MIO}^\star$ is an optimal solution of \eqref{eq:MIPII} and $\bm{\beta}_{persp}^\star$ is optimal to \eqref{eq:persp}.
\end{theorem}
{\color{black}We note that Problem \eqref{eq:MIPII} is a restatement of Problem \eqref{eqn:l0l2_training}.}

\proof{Proof of Theorem \ref{prop:kfoldspread}}
\rev{Let 
\begin{align}\label{persp:peturbed}
f(\bm{\beta})&:= \min_{\bm{z} \in [0, 1]^p: \bm{e}^\top \bm{z}\leq \tau}\quad \underbrace{\|\bm{X}\bm{\beta}-\bm{y}\|_2^2+\frac{\epsilon}{2}\|\bm{\beta}\|_2^2}_{=q(\bm{\beta})}+\underbrace{\frac{\gamma-\epsilon}{2}\sum_{i \in [p]}\frac{\beta_i^2}{z_i}}_{=r(\bm{\beta})}
\end{align}
denote the objective value of the perspective relaxation at a given $\bm{\beta}$, which can be decomposed in a quadratic part $q(\bm{\beta})$ and a convex nonlinear part $r(\bm{\beta})$. {\color{black}Since $q(\bm{\beta})$ can be rewritten as its second-order Taylor series expansion without loss of generality, we} find that for any $\bm{\beta}\in \mathbb{R}^p$,
\begin{align*}
q(\bm{\beta})&=q(\bm{\beta^\star_{\text{persp}}})+\nabla q(\bm{\beta^\star_{\text{persp}}})^\top (\bm{\beta}-\bm{\beta^\star_{\text{persp}}})+(\bm{\beta}-\bm{\beta^\star_{\text{persp}}})^\top\left(\bm{X^\top X}+\frac{\epsilon}{2}\mathbb{I}\right)(\bm{\beta}-\bm{\beta^\star_{\text{persp}}}) \text{, and}\\
r(\bm{\beta})&\geq r(\bm{\beta^\star_{\text{persp}}})+\bm{s}^\top (\bm{\beta}-\bm{\beta^\star_{\text{persp}}})
\end{align*}
for $\bm{s}\in \partial r(\bm{\beta^\star_{\text{persp}}})$ (the subdifferential of $r$). Moreover, since $\bm{\beta^\star_{\text{persp}}}$ is a minimizer of $f=q+r$, there exists $\bar{\bm{s}}\in \partial r(\bm{\beta^\star_{\text{persp}}})$ such that $\bar{\bm{s}}+\nabla q(\bm{\beta^\star_{\text{persp}}})=\bm{0}$. Adding the two inequalities, we find that
$$(\bm{\beta}-\bm{\beta^\star_{\text{persp}}})^\top\left(\bm{X^\top X}+\frac{\epsilon}{2}\mathbb{I}\right)(\bm{\beta}-\bm{\beta^\star_{\text{persp}}})\leq f(\bm{\beta})-f(\bm{\beta^\star_{\text{persp}}}).$$
Finally, setting $\bm{\beta}=\bm{\beta}_{MIO}^\star$ and using that $f(\bm{\beta}_{MIO}^\star)\leq \bar u$, we obtain the desired result.}
\Halmos
\endproof


\ag{Using Theorem~\ref{prop:kfoldspread}, we can compute bounds on {\color{black}$h_j(\gamma,\tau):=\sum_{i \in \mathcal{N}_j}\left(y_i -\bm{x_i}^\top \bm{\beta}^{(\mathcal{N}_j)}(\gamma, \tau)\right)^2$} in \eqref{eqn:crossValError} by solving problems of the form 
\begin{subequations}\label{eq:boundComputation}
\begin{align}
\min/\max_{\rev{\bm{\beta}\in \mathbb{R}^p}}\quad&\sum_{i\in \mathcal{N}_j} \left(y_i -\bm{x_i}^\top \bm{\beta}\right)^2 \label{eqn:9a}\\
\text{s.t.}\quad&(\bm{\beta_{\text{persp}}}^{(\mathcal{N}_j)}- \bm{\beta})^\top \left(\bm{X}^{(\mathcal{N}_j)^\top} \bm{X}^{(\mathcal{N}_j)}+\frac{\epsilon}{2}\mathbb{I}\right) (\bm{\beta_{\text{persp}}}^{(\mathcal{N}_j)}- \bm{\beta}) \leq (\bar u^{(\mathcal{N}_j)}-\zeta_{\text{persp}}^{(\mathcal{N}_j)}),\label{eqn:9b}
\end{align}
\end{subequations}
where $\bm{\beta_{\text{persp}}}^{(\mathcal{N}_j)}$ and $\zeta_{\text{persp}}^{(\mathcal{N}_j)}$ are the optimal solution and objective value of the perspective relaxation with fold $\mathcal{N}_j$ removed, and $\bar u^{(\mathcal{N}_j)}$ is an associated upper bound. {\color{black}Indeed, applying Theorem \ref{prop:kfoldspread} to the lower-level training problem with fold $\mathcal{N}_j$ removed implies that the true optimizer $\bm{\beta}^{(\mathcal{N}_j)}(\gamma,\tau)$ satisfies constraint \eqref{eqn:9b}, since $\bar u^{(N_j)}$ is an upper bound on the corresponding leave-fold-out MIO objective. Therefore, $h_j(\gamma,\tau)=\sum_{i\in N_j}(y_i-x_i^\top\beta^{(N_j)}(\gamma,\tau))^2$
is the value of the objective in \eqref{eqn:9a} at a feasible point of \eqref{eqn:9b}, so the minimum and maximum in \eqref{eqn:9a}--\eqref{eqn:9b} give valid lower and upper bounds on $h_j(\gamma,\tau)$, respectively.} Bounds for the function $h(\gamma,\tau)$ then immediately follow by simply adding the bounds associated with $h_j(\gamma,\tau)$ for all $j\in [k]$. 
}

\begin{remark}[Computability of \ag{the} bounds]\label{rem:computability}
	Observe that a lower bound on the $k$-fold error can easily be computed by solving a convex quadratically constrained quadratic problem, while an upper bound can be computed by noticing that the maximization problem \eqref{eq:boundComputation} is a trust region problem in $\bm{\beta}$, which can be reformulated as a semidefinite problem \citep{hazan2016linear}. One could further tighten these bounds by imposing a sparsity constraint on $\bm{\beta}$, but this may not be practically tractable.
\end{remark}



\subsection{Closed-form Bounds on the Prediction Spread}
\ag{While solving the perspective relaxation \eqref{eq:persp} is necessary to solve the MIO \eqref{eq:MIPII} via branch-and-bound (in particular, the perspective relaxation is the root node in a branch-and-bound scheme \citep{mazumder2023subset}), the additional two optimization problems \eqref{eq:boundComputation} are not. Moreover, solving trust-region problems can be expensive in large-scale problems. Accordingly, in this section, we present alternative bounds that may be weaker, but can be obtained in closed form. In numerical experiments (Section \ref{sec:numres}), these closed-form bounds {\color{black}typically} reduce the number of MIOs that need to be solved by {\color{black}50\%–80\%} when compared to grid search.}

\begin{theorem}\label{prop:l10UB}
	Given any vector $\bm{x}\in \mathbb{R}^p$ and any bound
	\begin{equation}\label{eq:MIPI}
			\bar u\geq \min_{\bm{\beta}\in \mathbb{R}^p}\; \|\bm{X}\bm{\beta}-\bm{y}\|_2^2+\frac{\gamma}{2}\|\bm{\beta}\|_2^2
			\text{\rm\ s.t. }\; \|\bm{\beta}\|_0\leq \tau,
	\end{equation} 
	the inequalities
	\small
	$$\bm{x^\top\beta}_{persp}^\star-\sqrt{\left(\bar u-\zeta_{\text{persp}}\right)\bm{x^\top}\left(\bm{X^\top X}+\frac{\epsilon}{2}\mathbb{I}\right)^{-1}\bm{x}}\leq \bm{x}^\top \bm{\beta}_{MIO}^\star\leq  \bm{x^\top\beta}_{persp}^\star+\sqrt{\left(\bar u-\zeta_{\text{persp}}\right)\bm{x^\top}\left(\bm{X^\top X}+\frac{\epsilon}{2}\mathbb{I}\right)^{-1}\bm{x}}$$\normalsize
	hold, where $\bm{\beta}_{MIO}^\star$ is an optimal solution of \eqref{eq:MIPI} and $\bm{\beta}_{persp}^\star$ is optimal to \eqref{eq:persp}.
\end{theorem}

\proof{Proof of Theorem~\ref{prop:l10UB}}
From Theorem~\ref{prop:kfoldspread}, we have the inequality
\begin{align}
(\bm{\beta^\star_{\text{persp}}}- \bm{\beta}^\star_{\text{MIO}})^\top \left(\bm{X}^\top \bm{X}+\frac{\epsilon}{2}\mathbb{I}\right) (\bm{\beta^\star_{\text{persp}}}- \bm{\beta}^\star_{\text{MIO}}) \leq (\bar{u}-\zeta_{\text{persp}}).
	\end{align}
 By the Schur Complement Lemma \citep[see, e.g.,][]{boyd1994linear}, this is equivalent to
 \begin{align*}
     (\bar{u}-\zeta_{\text{persp}}) \left(\bm{X}^\top \bm{X}+\frac{\epsilon}{2}\mathbb{I}\right)^{-1}\succeq (\bm{\beta^\star_{\text{persp}}}- \bm{\beta}^\star_{\text{MIO}})(\bm{\beta^\star_{\text{persp}}}- \bm{\beta}^\star_{\text{MIO}})^\top.
 \end{align*}
 Next, we can left/right multiply this expression by an arbitrary matrix $\bm{W} \in \mathbb{R}^{m \times p}$. This gives:
 \begin{align*}
     (\bar{u}-\zeta_{\text{persp}}) \bm{W}\left(\bm{X}^\top \bm{X} +\frac{\epsilon}{2}\mathbb{I}\right)^{-1}\bm{W}^\top\succeq (\bm{W}\bm{\beta^\star_{\text{persp}}}- \bm{W}\bm{\beta}^\star_{\text{MIO}})(\bm{W}\bm{\beta^\star_{\text{persp}}}- \bm{W}\bm{\beta}^\star_{\text{MIO}})^\top.
 \end{align*}
In particular, setting $\bm{W}=\bm{x}^\top$ for a vector $\bm{x} \in \mathbb{R}^{p}$ gives the inequality
 \begin{align*}
     (\bar{u}-\zeta_{\text{persp}}) \bm{x}^\top \left(\bm{X}^\top \bm{X}+\frac{\epsilon}{2}\mathbb{I}\right)^{-1} \bm{x}\geq (\bm{x}^\top(\bm{\beta^\star_{\text{persp}}}- \bm{\beta}^\star_{\text{MIO}}))^2,
 \end{align*}
 which we rearrange to obtain the result.
\Halmos
	
\endproof

\begin{corollary}
    For any $\bm{W} \in \mathbb{R}^{m \times p}$ we have that
 \begin{align*}
     (\bar{u}-\zeta_{\text{persp}}) \mathrm{tr}\left(\bm{W} \left(\bm{X}^\top \bm{X}+\frac{\epsilon}{2} \bm{\mathbb{I}}\right)^{-1} \bm{W}^\top\right)\geq \Vert 
     \bm{W}(\bm{\beta^\star_{\text{persp}}}- \bm{\beta}^\star_{\text{MIO}})\Vert_2^2\rev{.}
 \end{align*}
\end{corollary}

Applying Theorem~\ref{prop:l10UB} to the problem 
\begin{equation}\label{eq:MIOFoldJ}	\bar u^{(\mathcal{N}_j)}\geq \min_{\bm{\beta}\in \mathbb{R}^p}\; \|\bm{X}^{(\mathcal{N}_j)}\bm{\beta}-\bm{y}^{(\mathcal{N}_j)}\|_2^2+\frac{\gamma}{2}\|\bm{\beta}\|_2^2
\text{ s.t. }\; \|\bm{\beta}\|_0\leq \tau,\end{equation} \rev{where $\bm{x}=\bm{x_i}$ was chosen as the reference vector,} we have the bounds
\begin{align}\label{eq:foldij}
	\underline{\xi}_{i,j}:=\bm{x_i^\top}\bm{\beta}_{persp}^{{\color{black}(\mathcal{N}_j)}\star}-\sqrt{\bm{x_i^\top}\left(\bm{X}^{(\mathcal{N}_j)^\top} \bm{X}^{({\color{black}\mathcal{N}_j})}+\frac{\epsilon}{2}\mathbb{I}\right)^{-1}\bm{x_i}\left(\bar u^{(\mathcal{N}_j)}-\zeta^{(\mathcal{N}_j)}\right)},\\
	\bar{\xi}_{i,j}:=\bm{x_i^\top}\bm{\beta}_{persp}^{{\color{black}(\mathcal{N}_j)}\star}+\sqrt{\bm{x_i^\top}\left(\bm{X}^{(\mathcal{N}_j)^\top} \bm{X}^{(\mathcal{N}_j)}+\frac{\epsilon}{2}\mathbb{I}\right)^{-1}\bm{x_i}\left(\bar u^{(\mathcal{N}_j)}-\zeta^{(\mathcal{N}_j)}\right)}
\end{align}
where $\underline{\xi}\leq \bm{x_i^\top}\bm{\beta}_{MIO}^{{\color{black}(\mathcal{N}_j)}\star}\leq \bar\xi$\rev{, where $\bm{\beta}_{MIO}^{{\color{black}(\mathcal{N}_j)}\star}$ is the optimal solution of \eqref{eq:MIOFoldJ}.} \rev{We can then compute bounds on the $i$th prediction error associated with fold $j$, namely $$\nu_{i,j}=(\bm{x_i^\top}\bm{\beta}_{MIO}^{{\color{black}(\mathcal{N}_j)}\star}-y_i)^2,$$ as formalized in Corollary~\ref{corr:perspUB}.}

\begin{corollary}\label{corr:perspUB}
	We have the following bounds on the $i$th prediction error \ag{associated with fold $j$}
	\begin{align}\label{eqn:perspbounds}
		\max\left((y_i-\underline{\xi}_{i,j})^2, (y_i-\bar{\xi}_{i,j})^2\right) \geq \ag{\nu_{i,j}(\gamma, \tau)} \geq \begin{cases} (y_i-\underline{\xi}_{i,j})^2 \quad & \text{if} \ y_i < \underline{\xi}_{i,j}\\
			0 \quad & \text{if} \ y_i \in [\underline{\xi}_{i,j}, \bar{\xi}_{i,j}],\\
			(\bar{\xi}_{i,j}-y_i)^2 \quad & \text{if} \ y_i > \bar{\xi}_{i,j}.
		\end{cases}
	\end{align}
\end{corollary}
Moreover, since $h(\gamma,\tau)=\frac{1}{n}\sum_{j=1}^k\sum_{i\in \mathcal{N}_j}\nu_{i,j}(\gamma,\tau)$, we can compute lower and upper bounds on the $k$-fold cross-validation error \ag{by adding the individual bounds. Observe that the bounds computed by summing disaggregated bounds could be substantially worse than those obtained by letting $\bm{W}$ be a matrix with all omitted {\color{black}rows} in the $j$th fold of $\bm{X}$ in the proof of Theorem \ref{prop:l10UB}. {\color{black}Nonetheless, the approach outlined here may be the only feasible one in large-scale instances, as it is obtained directly from the perspective relaxation without solving additional optimization problems, whereas an aggregated approach would involve solving an auxiliary semidefinite optimization problem.} Despite the loss in quality, we show in our computational sections that (combined with the methods discussed in \S\ref{sec:corddescent}), the disaggregated bounds are sufficient to lead to a 50\%–80\% reduction in the number of MIO{\color{black}s} solved with respect to grid search.}

We conclude this subsection with \rev{three} remarks.

\rev{\begin{remark}[Choice of the strong convexity parameter]While it may appear that large values $\epsilon$ result in tighter bounds in Theorem~\ref{prop:l10UB}, we point out that large values also negatively affect the quality of the lower bound $\zeta_{\text{persp}}$. Note that if $\bm{X}^\top\bm{X}$ is invertible, then we can also set $\epsilon=0$ -- the approach we use if $n>p$. 
\end{remark}}

\begin{remark}[Relaxation Tightness]
	If the perspective relaxation is tight, as occurs 
 when $n$ is sufficiently large under certain assumptions on the data generation process \citep{pilanci2015sparse, reeves2019all}
 then 	$\underline{\xi}=\bar{\xi}=\bm{x_i^\top}\bm{\beta}_{persp}^\star$, and Corollary \ref{corr:perspUB}'s bounds on the cross-validation error are tight {\color{black}by definition}. Otherwise, as pointed out in Remark~\ref{rem:intuition}, \eqref{eqn:perspbounds}'s bound quality depends on the tightness of the relaxation and on how close the features $\bm{x}_i$ are to the rest of the data. 
\end{remark}

\begin{remark}[Intuition]\label{rem:intuition}
	Theorem~\ref{prop:l10UB} states that $\bm{x}^\top \bm{\beta}_{MIO}^\star\approx\bm{x^\top\beta}_{persp}^\star$, where the approximation error is determined by two components. The quantity $\sqrt{\bar u-\zeta_{\text{persp}}}$ is related to the strength of the perspective relaxation, with a stronger relaxation resulting in a better approximation. The quantity $\sqrt{\bm{x^\top}\left(\bm{X^\top X}+\frac{\epsilon}{2}\mathbb{I}\right)^{-1}\bm{x}}$ is related to the likelihood that $\bm{x}$ is generated from the same distribution as the rows of $\bm{X}$, with larger likelihoods resulting in better approximations. Indeed, if $n>p$, each column of $\bm{X}$ has $0$ mean but has not been standardized, and each row of $\bm{X}$ is generated iid from a multivariate Gaussian distribution, then $\frac{n(n-1)}{n+1}\bm{x^\top}\left(\bm{X^\top X}\right)^{-1}\bm{x}\sim T^2(p, n-1)$ is Hotelling's two-sample T-square test statistic \citep{hotelling1931generalization}, used to test whether $\bm{x}$ is generated from the same Gaussian distribution. Note that if $\bm{x}$ is drawn from the same distribution as the rows of $\bm{X}$ (as may be the case in cross-validation), then $\mathbb{E}\left[\bm{x^\top}\left(\bm{X^\top X}\right)^{-1}\bm{x}\right]=\frac{p(n+1)}{n(n-p-2)}$ {\color{black}for $n>p+2$}.
\end{remark}

\subsection{\ag{Further Improvements for Lower Bounds}}\label{sec:genBounds}

Corollary~\ref{corr:perspUB} implies we obtain valid upper and lower bounds on {\color{black}the $k$-fold cross-validation loss }$h$ at a given hyperparameter combination $\gamma, \tau$ after solving $k$ perspective relaxations and computing $n$ terms of the form $$\sqrt{\bm{x_i^\top}\ag{\left(\bm{X}^{(\mathcal{N}_j)^\top} \bm{X}^{(\mathcal{N}_j)}+\frac{\epsilon}{2}\mathbb{I}\right)^{-1}}\bm{x_i}}.$$ 

A drawback of Corollary~\ref{corr:perspUB} is that if $\bm{x_i^\top}\bm{\beta}_{persp}^\star\approx y_i$ for each $i \in \mathcal{N}_j$, i.e., the prediction of the perspective relaxation (without \ag{the $j$th fold}) is close to the response associated with point $i$, then Corollary~\ref{corr:perspUB}'s lower bound is $0$. \ag{A similar situation can happen with the stronger bounds for $h_j(\gamma,\tau)$ obtained from Theorem~\ref{prop:kfoldspread} and Problem \eqref{eq:boundComputation}.} We now propose a different bound on \ag{$h_j(\gamma,\tau)$}, which is sometimes effective in this circumstance.

First, define the function ${\color{black}F}(\gamma,\tau)$ to be the in-sample training error \ag{without removing any folds and }with parameters $(\gamma, \tau)$,
\begin{equation*}
{\color{black}F}(\gamma, \tau) := \frac{1}{n} \sum_{i=1}^n (y_i -\bm{x_i}^\top \bm{\beta}(\gamma, \tau))^2\quad
\text{s.t.} \quad  \bm{\beta}(\gamma, \tau) \in \argmin_{\bm{\beta} \in \mathbb{R}^p: \ \Vert \bm{\beta}\Vert_0 \leq \tau} \ \frac{\gamma}{2}\Vert \bm{\beta}\Vert_2^2 +\Vert \bm{X}\bm{\beta}-\bm{y}\Vert_2^2,
\end{equation*}
and let \ag{${\color{black}F}_{\rev{j}}(\gamma, \tau):=\sum_{i\in \mathcal{N}_j}\left(y_i -\bm{x_i}^\top \bm{\beta}(\gamma,\tau)\right)^2$} denote the training error \ag{associated with the $j$th fold}, with $\ag{1/n \sum_{j=1}^k {\color{black}F}_{\rev{j}}(\gamma, \tau)}={\color{black}F}(\gamma,\tau)$. Observe that evaluating \rev{$h$} involves solving \ag{$k$} MIOs, while evaluating ${\color{black}F}$ requires solving one. 

\begin{proposition}\label{prop:LB}
For any $\gamma{\color{black}>} 0$, any $\tau \in[p]$ \ag{and any $j\in [k]$}, ${\color{black}F}_j(\gamma,\tau)\leq h_j(\gamma,\tau)$. Moreover, we have that ${\color{black}F}(\gamma,\tau) \leq h(\gamma, \tau)$.
\end{proposition}
\proof{Proof of Proposition \ref{prop:LB}}
Given ${\color{black}j}\in [k]$, consider the following two optimization problems
	\begin{align}
		\min_{\bm{\beta}\in \mathbb{R}^p: \Vert \bm{\beta}\Vert_0 \leq \tau}\sum_{i=1}^n(y_i -\bm{x_i}^\top \bm{\beta})^2+\frac{\gamma}{2}\Vert \bm{\beta}\Vert_2^2\label{eq:prob_full}\\
	\min_{\bm{\beta}\in \mathbb{R}^p: \Vert \bm{\beta}\Vert_0 \leq \tau}\sum_{\ag{i \not\in \mathcal{N}_j}}(y_i -\bm{x_i}^\top \bm{\beta})^2+\frac{\gamma}{2}\Vert \bm{\beta}\Vert_2^2,\label{eq:prob_noK}
	\end{align}
let $\bm{\beta^\star}$ be an optimal solution of \eqref{eq:prob_full}, and let $\bm{\beta^{\color{black}j}}$ be an optimal solution of \eqref{eq:prob_noK}. Since \begin{align*}&\sum_{\ag{i \not\in \mathcal{N}_j}}(y_i -\bm{x_i}^\top \bm{\beta^{\color{black}j}})^2+\frac{\gamma}{2}\Vert \bm{\beta^{\color{black}j}}\Vert_2^2\leq \sum_{\ag{i \not\in \mathcal{N}_j}}(y_i -\bm{x}^\top \bm{\beta^\star})^2+\frac{\gamma}{2}\Vert \bm{\beta^\star}\Vert_2^2,\quad \text{ and }\\
&\sum_{\ag{i \not\in \mathcal{N}_j}}(y_i -\bm{x}^\top \bm{\beta^{\color{black}j}})^2+\ag{\sum_{i\in \mathcal{N}_j}}(y_{\color{black}i} -\bm{x_{\color{black}i}}^\top \bm{\beta^{\color{black}j}})^2+\frac{\gamma}{2}\Vert \bm{\beta^{\color{black}j}}\Vert_2^2\geq \sum_{\ag{i \not\in \mathcal{N}_j}}(y_i -\bm{x_i}^\top \bm{\beta^\star})^2+\ag{\sum_{i\in \mathcal{N}_j}}(y_{\color{black}i} -\bm{x_{\color{black}i}}^\top \bm{\beta^\star})^2+\frac{\gamma}{2}\Vert \bm{\beta^\star}\Vert_2^2,
\end{align*}
we conclude that $\ag{\sum_{i\in \mathcal{N}_j}}(y_{\color{black}i} -\bm{x_{\color{black}i}}^\top \bm{\beta^\star})^2\leq  \ag{\sum_{i\in \mathcal{N}_j}}(y_{\color{black}i} -\bm{x_{\color{black}i}}^\top \bm{\beta^{\color{black}j}})^2$. The result immediately follows. \Halmos
\endproof


Next, we develop a stronger bound on the $k$-fold error, by observing that our original proof technique relies on interpreting the optimal solution when training on the entire dataset as a feasible solution when leaving out the \ag{$j$}th fold, and that this feasible solution can be improved to obtain a tighter lower bound. Therefore, {given \ag{any} $\bm{z}\in \{0,1\}^p$,} let us define the function:{\color{black}
\begin{align*}
    f^{(\ag{\mathcal{N}_j})}(\bm{z}):=\min_{\bm{\beta} \in \mathbb{R}^p} \quad \frac{\gamma}{2}\sum_{j \in [p]}\beta_j^2+\Vert \bm{X}^{(\ag{\mathcal{N}_j})}\bm{\beta}-\bm{y}^{(\ag{\mathcal{N}_j})}\Vert_2^2 \quad \text{s.t.} \quad \beta_j=0 \ \text{if} \ z_j=0 \ \forall j \in [p],
\end{align*}}
to be the optimal training loss (including regularization) when we leave out the \ag{$j$th fold} and have the binary support vector $\bm{z}$. Then, fixing $\gamma, \tau$ and letting $u^\star$ denote the optimal objective value of \eqref{eq:prob_full}, i.e., the optimal training loss on the entire dataset (including regularization) and $\bm{\beta}^{(\ag{\mathcal{N}_j})}(\bm{z})$ denote an optimal choice of $\bm{\beta}$ for this $\bm{z}$, we have the following result:
\begin{proposition}\label{prop2}
For any $\tau$-sparse binary vector $\bm{z}$, the following inequality holds:
\begin{align}
    u^\star\leq f^{(\ag{\mathcal{N}_j})}(\bm{z})+\ag{\sum_{i\in \mathcal{N}_j}}\left(y_i-\bm{x}_i^\top \bm{\beta}^{(\ag{\mathcal{N}_j})}(\bm{z})\right)^2
\end{align}
\end{proposition}
\proof{Proof of Proposition \ref{prop2}}
The right-hand side of this inequality corresponds to the objective value of a feasible solution to \eqref{eq:prob_full}, while {\color{black}$u^\star$} is the optimal objective value of \eqref{eq:prob_full}. \Halmos
\endproof
\begin{corollary}\label{corr:lowerbound2}
Let $\bm{z}$ denote a $\tau$-sparse binary vector. Then, we have the following bound on the $j$th partial \ag{cross-validation} error:
\begin{align}\label{eq:LOOLower}
    h_j(\gamma, \tau) \geq u^\star-f^{(\ag{\mathcal{N}_j})}(\bm{z}).
\end{align}
\end{corollary}
\proof{Proof of Corollary \ref{corr:lowerbound2}}
The right-hand side of this bound is maximized by setting $\bm{z}$ to be a binary vector which minimizes $f^{(\ag{\mathcal{N}_j})}(\bm{z})$, and therefore this bound is valid for any $\bm{z}$. \Halmos
\endproof
We close this section with two remarks:
\begin{remark}[Bound quality]
{\color{black}Bound \eqref{eq:LOOLower} is designed to complement, rather than dominate, the lower bound in \eqref{eqn:perspbounds}}. Indeed, bound~\eqref{eq:LOOLower} is at least as strong as ${\color{black}F}_j(\gamma, \tau)$ with $\bm{z}$ encoding an optimal choice of support in \eqref{eq:prob_full}: if $\bm{\beta}^{(\ag{\mathcal{N}_j})}(\bm{z})$ solves \eqref{eq:prob_full}, then both bounds agree and equal $h_j(\gamma, \tau)$ but otherwise \eqref{eq:LOOLower} is strictly stronger. Moreover, since ${\color{black}F}_j(\gamma, \tau)$ is typically nonzero, then the bound \eqref{eq:LOOLower} is positive as well and can improve upon the lower bound in \eqref{eqn:perspbounds}. {\color{black}However,} it is easy to construct examples {\color{black}in which} the lower bound in \eqref{eqn:perspbounds} is stronger than \eqref{eq:LOOLower} and vice versa, {\color{black}so} neither lower bound dominates the other{\color{black}; see Section \ref{ec.nondominanceofrelaxations}}.
\end{remark}
\begin{remark}[Computational efficiency]
	Computing lower bound \eqref{eq:LOOLower} for each \ag{$j\in [k]$} requires solving at least one MIO, corresponding to \eqref{eq:prob_full}, which is a substantial improvement over the \ag{$k$} MIOs required to compute $h$ but may still be an expensive computation. However, using any lower bound on $u^\star$, for example, corresponding to the optimal solution of a perspective relaxation, gives valid lower bounds. Therefore, in practice, we suggest using a heuristic instead to bound $h_j$ from below, e.g., rounding a perspective relaxation {\color{black}as suggested by \cite{xie2020scalable,bertsimas2021unified}, building upon the work of \cite{pilanci2015sparse}}.
\end{remark}


\section{Optimizing the Cross-Validation Loss}\label{sec:corddescent}
In this section, we present an efficient \rev{alternating minimization} scheme that identifies (approximately) optimal hyperparameters $(\gamma, \tau)$ with respect to the {\color{black}$k$-fold cross-validation error as previously defined in \eqref{prob:upperlevelofv_validation}, \eqref{eqn:crossValError}:} 
\begin{equation}\label{prob:cyclic1}\color{black}
	h(\gamma, \tau):= \frac{1}{n} \sum_{j \in [k]} h_j(\gamma, \tau)
\end{equation}
by iteratively minimizing $\tau$ and $\gamma$. \rev{Specifically}, with initialization $\tau_0, \gamma_0$, we repeatedly solve the following two optimization problems:
\begin{align}\label{prob:cyclic2a}
    \tau_t \in \argmin_{\tau \in[p]} \quad h(\gamma_t, \tau),\\
    \gamma_{t+1} \in \arg \min_{\gamma >0 } \quad h(\gamma, \tau_t),\label{prob:cyclic2b}
\end{align}
until we either detect a cycle, converge to a locally optimal solution {\color{black} to \eqref{prob:cyclic1}, or exceed a user-imposed limit on the number of iterations for this alternating minimization procedure}. To develop this scheme, in Section \ref{ssec:parametric} we propose an efficient technique for solving Problem {\color{black}\eqref{prob:cyclic2a}} {(\color{black}Algorithm \ref{alg:parametricK})}, 
and in Section \ref{ssec:parametric2} we propose an efficient technique for (approximately) solving Problem {\color{black}\eqref{prob:cyclic2b}}. {\color{black}Our overall scheme alternates between solving the two minimization problems.} Accordingly, our scheme could also be used to identify a {\color{black}locally} optimal choice of $\gamma$ if $\tau$ is already known, e.g., in a context where regulatory constraints specify the number of features that may be included in a model. 

Our overall approach is motivated by three key observations. First, we design a method that obtains local, rather than global, minima, because $h$ is a highly non-convex function and even evaluating $h$ requires solving {\color{black}$k$} MIOs, which suggests that global minima of $h$ may not be attainable in a practical amount of time at scale. Second, we use \rev{alternating minimization} to seek local minima because if either $\tau$ or $\gamma$ is fixed, it is possible to efficiently optimize the remaining hyperparameter with respect to $h$ by leveraging the convex relaxations developed in the previous section. Third, we should expect our \rev{alternating minimization} scheme to perform well in practice, because similar schemes are highly effective in other machine learning contexts, e.g., solving certain matrix completion problems in polynomial time \citep{mazumder2011sparsenet, cifuentes2022polynomial}.


\subsection{Parametric Optimization of $k$-fold \rev{Error} With Respect to Sparsity}\label{ssec:parametric}

Consider the following optimization problem, where $\gamma$ is fixed here and throughout this subsection:
{\color{black}
\begin{align}\label{eq:paramK}
\min_{\tau\in [p]} \quad h({\gamma}, \tau):=& {\color{black}\frac{1}{n}}\sum_{j \in [k]}\sum_{i \in \mathcal{N}_j} (y_i-\bm{x}_i^\top \bm{\beta}^{(\mathcal{N}_j)})^2,\\
\text{s.t.} \quad  & \bm{\beta}^{(\mathcal{N}_j)} \in \argmin_{\bm{\beta} \in \mathbb{R}^p: \ \Vert\bm{\beta}\Vert_0\leq \tau} \ \frac{\gamma}{2}\Vert \bm{\beta}\Vert_2^2 +\Vert \bm{X}^{(\mathcal{N}_j)}\bm{\beta}-\bm{y}^{(\mathcal{N}_j)}\Vert_2^2\quad \forall \color{black}j \in [k].\nonumber
\end{align}}
{\color{black}Note that solving Problem \eqref{eq:paramK} corresponds to performing one iteration of the alternating minimization scheme described at the start of the section, namely minimizing $h$ with respect to $\tau$ with $\gamma$ fixed.

Problem \eqref{eq:paramK}} can be solved by complete enumeration, i.e., for each $\tau \in[p]$, we compute an optimal $\beta^{(\mathcal{N}_j)}$ for each $j\in [k]$ by solving an MIO. This involves solving $kp$ MIOs, which is extremely expensive at scale. We now propose a technique for minimizing $h$ without solving all these MIOs: 

{Algorithm \ref{alg:parametricK} has two main phases, {\color{black}each implemented as a loop over sparsity budgets $\tau$ and then over folds $j$}. In the first phase, we construct valid lower and upper bounds on $h_{\rev{j}}(\rev{\gamma,}\tau)$ for each $\rev{j\in [k]}$ and each $\tau$ without solving any MIOs. We begin by solving, for each potential sparsity budget $\tau\in [p]$, the perspective relaxation with all data points included. Call this relaxation's objective value $\bar v_{\tau}$. We then solve each perspective relaxation that arises after omitting one data fold $\mathcal{N}_j: j \in [k]$, with objective values $v_{\tau,\rev{j}}$ and solutions $\bm{\beta}_{\tau,\rev{j}}$. 
Next, we compute lower and upper bounds on the $k$-fold error $h_{\rev{j}}(\rev{\gamma,}\tau)$ using the methods derived in Section~\ref{sec:lowerbounds}, which are summarized in the routine \texttt{compute\_bounds} described in Algorithm~\ref{alg:bounds}. By solving $\mathcal{O}(pk)$ relaxations (and no MIOs), we have upper and lower estimates on the $k$-fold error that are often accurate in practice, as described by Theorem~\ref{prop:l10UB}. 

After completing the first loop in Algorithm~\ref{alg:parametricK}, one may already terminate the algorithm. Indeed, according to our numerical experiments in Section \ref{sec:numres}, this already provides high-quality solutions. Alternatively, one may proceed with the second phase of Algorithm~\ref{alg:parametricK} and solve \eqref{prob:cyclic2a} to optimality, at the expense of solving (a potentially large number of) MIOs.

In the second phase, Algorithm~\ref{alg:parametricK} identifies the cardinality $\tau^\star$ with the best lower bound (and thus, in an optimistic scenario, the best potential value). Then, it identifies the \rev{index $j^*\in[k]$ such that partition $\mathcal{N}_{j^\star}$ has} the largest uncertainty around the $k$-fold estimate $h_{\rev{j^\star}}(\rev{\gamma,}\tau^\star)$, and solves an MIO to compute the exact partial $k$-fold error. 
This process is repeated until \eqref{eq:paramK} is solved \rev{within a certain prescribed optimality tolerance {\color{black}$\text{opt}_{\rm tol}$}, or a suitable termination condition (e.g., a limit on computational time) 
is met. {\color{black}Note that if an MIO were solved for all $\tau$ and all folds of the data, then we could immediately solve \eqref{eq:paramK}, and thus Algorithm~\ref{alg:parametricK} terminates in $pk$ iterations of the second phase in the worst case.}

}


\normalem

\RestyleAlgo{ruled}
\begin{algorithm}[h!] 
	\caption{Computing optimal sparsity parameter for $k$-fold error}
	\label{alg:parametricK}
	\KwData{$\gamma$: $\ell_2^2$ regularization parameter; $\text{opt}_{\text{tol}}>0$: desired optimality tolerance; $r$: budget on number of MIOs}
	\KwResult{Cardinality with best estimated $k$-fold error} 
        \For{$\tau \in[p]$}{
        $\bar v_{\tau} \gets \min_{\bm{\beta}\in \mathbb{R}^p,\bm{z}\in [0,1]^p}\;\|\bm{X\beta}-\bm{y}\|_2^2+\frac{\gamma}{2}\sum_{i=1}^p\beta_i^2/z_i \text{  s.t.  }\bm{e}^\top \bm{z}\leq \tau$\\
			\For{$j\in [k]$}{
                $v_{\tau,\rev{j}} \gets\min_{\bm{\beta}\in \mathbb{R}^p,\bm{z}\in [0,1]^p}\;\|\bm{X}^{(\mathcal{N}_j)}\bm{\beta}-\bm{y}^{(\mathcal{N}_j)}\|_2^2+\frac{\gamma}{2}\sum_{i=1}^p\beta_i^2/z_i \text{  s.t.  }\bm{e}^\top \bm{z}\leq \tau$\\
                $\bm{\beta}_{\tau,\rev{j}}, \bm{z}_{\tau,j} \in \argmin_{\bm{\beta}\in \mathbb{R}^p,\bm{z}\in [0,1]^p}\;\|\bm{X}^{(\mathcal{N}_j)}\bm{\beta}-\bm{y}^{(\mathcal{N}_j)}\|_2^2+\frac{\gamma}{2}\sum_{i=1}^p\beta_i^2/z_i \text{  s.t.  }\bm{e}^\top \bm{z}\leq \tau$\\
                {\color{black}$\hat{\bm{z}}_{\tau, j} \gets \texttt{round}(\bm{z}_{\tau, j})$\tcp*{\color{black}set $\tau$ largest entries of $\bm{z}_{\tau, j}$ to $1$, remaining entries to $0$}
                $u_{\tau,\rev{j}} \gets\min_{\bm{\beta}\in \mathbb{R}^p}\;\|\bm{X}^{(\mathcal{N}_j)}\bm{\beta}-\bm{y}^{(\mathcal{N}_j)}\|_2^2+\frac{\gamma}{2}\sum_{i=1}^p\beta_i^2 \ \text{s.t.} \ \beta_i=0 \ \text{if} \ \hat{z}_i=0 \ \forall i \in [p]$ \\}
                $\zeta_{\rev{j}}^L(\tau),\zeta_{\rev{j}}^U(\tau)\gets \texttt{compute\_bounds}(\mathcal{N}_j,\bm{\beta}_{\tau,\rev{j}},\bar v_{\tau},v_{\tau,\rev{j}},u_{\tau,\rev{j}})$
		}\vspace{0.5mm}
	}
$LB\gets \min_{\tau \in[p]}\sum_{j \in [k]} \zeta_{\rev{j}}^L(\tau)$; 
$UB\gets \min_{\tau \in[p]}\sum_{j \in [k]} \zeta_{\rev{j}}^U(\tau)$ \tcp*{Bounds on $k$-fold} 
        $num\_mip\gets 0$\\
	\Repeat{$(UB-LB)/UB\rev{\leq} \text{opt}_{\text{tol}}$ or $num\_mip>r$}{
		$\rev{\bar\tau}\gets{\arg \min}_{\tau\in [p]} \sum_{\color{black}j \in [k]} \zeta_{\rev{j}}^L(\tau)$ \tcp*{Cardinality with best bound}
		$\rev{j^\star}\gets \arg \max_{j\in [k]} \{\zeta_{\rev{j}}^U(\rev{\bar\tau})-\zeta_{\rev{j}}^L(\rev{\bar\tau})\}$ \tcp*{Fold with largest $k$-fold uncertainty}{\color{black}
        $\bm{\beta}^\star \gets 
\argmin_{\bm{\beta}\in \mathbb{R}^p,\bm{z}\in \{0,1\}^p}\;\|\bm{X}^{(\mathcal{N}_j^*)}\bm{\beta}-\bm{y}^{(\mathcal{N}_j^*)}\|_2^2+\frac{\gamma}{2}\sum_{i=1}^p\beta_i^2/z_i \text{  s.t.  }\bm{e}^\top \bm{z}\leq \bar{\tau}$
\tcp*{Solve the training problem with fold $\mathcal{N}_{j^\star}$ left out.}
$h_{j^\star}(\gamma,\bar\tau) \gets \sum_{i \in \mathcal{N}_{j^\star}} \left(y_i - x_i^\top \beta^\star\right)^2$
\tcp*{Evaluate valid loss on held-out fold $\mathcal{N}_{j^\star}$.}
            $\zeta_{\rev{j^\star}}^L(\rev{\bar\tau})\gets h_{\rev{j^\star}}(\rev{\gamma,\bar\tau}),\; \zeta_{\rev{j^\star}}^U(\rev{\bar\tau})\gets h_{\rev{j^\star}}(\rev{\gamma,\bar\tau})$\\}	
            $LB\gets \min_{\tau \in[p]}\sum_{j \in [k]} \zeta_{\rev{j}}^L(\tau)$\\
            $UB\gets \min_{\tau \in[p]}\sum_{j \in [k]} \zeta_{\rev{j}}^U(\tau)$\\
            $num\_mip\gets num\_mip+1$\\[0.5em]
	}
        \Return{${\arg \min}_{\tau\in [p]} \sum_{j \in [k]} h_{\rev{j}}(\rev{\gamma},\tau)$}\tcp*{Cardinality with best error}
\end{algorithm}

\begin{algorithm}[h!] 
	\caption{$\texttt{compute\_bounds}(\mathcal{N}_j,\bm{\beta},\bar v,v,u)$}
	\label{alg:bounds}
	\KwData{$\mathcal{N}_j$: fold left out; $\bm{\beta}$: optimal solution of perspective relaxation with $\mathcal{N}_j$ left out; $\bar v$: lower bound of \rev{the objective value} of MIO with all data; $v$: optimal \rev{objective value} of perspective relaxation with $\mathcal{N}_j$ left out; $u$: upper bound of \rev{the objective value} of MIO with $\mathcal{N}_j$ left out{\color{black}; $\epsilon$: strong convexity parameter} }
	\KwResult{Lower and upper bounds on the $k$-fold error attributable to {\color{black}fold $j$}}
\For{\color{black}$i \in \mathcal{N}_j$}{
 $\underline{\xi_i}\gets\bm{x_i^\top\beta}-\sqrt{\bm{x_i^\top}\left(\bm{X}^{(\mathcal{N}_j)^\top} \bm{X}^{(\mathcal{N}_j)}+\frac{\epsilon}{2}\mathbb{I}\right)^{-1}\bm{x_i}\left(u-v\right)}$\\
$\overline{\xi_i}\gets\bm{x_i^\top\beta}+\sqrt{\bm{x_i^\top}\left(\bm{X}^{(\mathcal{N}_j)^\top} \bm{X}^{(\mathcal{N}_j)}+\frac{\epsilon}{2}\mathbb{I}\right)^{-1}\bm{x_i}\left(u-v\right)}$\\
  $\zeta_{\color{black}i}^L\gets {\color{black}0},\; \zeta_i^U\gets \max\{(y_i-\underline{\xi})^2, (\overline{\xi}-y_i)^2\}$\\
  \If{$\underline{\xi_i}> y_i$}{
   $\zeta_{\color{black}i}^L\gets \max\{\zeta_{\color{black}i}^L, (\underline{\xi}-y_i)^2\}$
  }
  \If{$\overline{\xi}< y_i$}{
   $\zeta_{\color{black}i}^L\gets \max\{\zeta_{\color{black}i}^L, (y_i-\overline{\xi})^2\}$
  }
}
	\Return{\color{black}$\left(\max(\bar v-u, \sum_{\color{black}i \in \mathcal{N}_j}\zeta_i^L), \sum_{\color{black}i \in \mathcal{N}_j}\zeta_i^U\right)$}
\end{algorithm}

To solve each MIO in Algorithm \ref{alg:parametricK}, we invoke a Generalized Benders Decomposition scheme \citep{geoffrion1972generalized}, which was specialized to sparse regression problems by \cite{bertsimas2020sparse}, enhanced with some ideas from the optimization literature summarized in the works \cite{bertsimas2020sparse2, hazimeh2020fast}. For the sake of conciseness, we defer these implementation details to Appendix \ref{sec:append.implementationdetails}.


\paragraph{Algorithm \ref{alg:parametricK} in Action: } Figure \ref{fig:algvarykinaction} depicts visually the lower and upper bounds on $\rev{h}$ from Algorithm \ref{alg:bounds} (left) and after running Algorithm \ref{alg:parametricK} to completion (right) on a synthetic sparse regression instance generated in the fashion {\color{black}described by \cite{bertsimas2020sparse2} and restated in Section \ref{append.syntheticdata} for completeness}, with $k=n$, $n=200, p=20$, $\gamma=1/\sqrt{n}$, $\tau_{\text{true}}=10$, $\rho=0.7$, $\nu=1$, where $\tau \in\{2, \ldots, 19\}${\color{black}, we have the tolerance parameters $r=kp$ and $\text{opt}_{\text{tol}}=10^{-2}$}, and using the outer-approximation method of \cite{bertsimas2020sparse} as our solver for each MIO with a time limit of $60$s. We observe that Algorithm \ref{alg:parametricK} solved $1694$ MIOs to identify the optimal $\tau$, which is a 53\% improvement \rev{over} complete enumeration. Interestingly, when $\tau=19$, the perspective relaxation is tight after omitting any fold of the data and we have tight bounds on the LOOCV error without solving any MIOs. In Section \ref{sec:compExactLOOCV}, we test Algorithm~\ref{alg:parametricK} on real datasets and find that it reduces the number of MIOs that need to be solved by 50-80\% with respect to complete enumeration. For more information on how the bounds evolve over time, we provide a GIF with one frame each time a MIO is solved at the link \url{https://drive.google.com/file/d/1EZdNwlV9sEEnludGGM7v2nGpB7tzZvz4/view?usp=sharing}.


\begin{figure}[h!]\centering
\begin{tikzpicture}
  \node[inner sep=0pt] (left) 
    {\includegraphics[width=.42\linewidth]{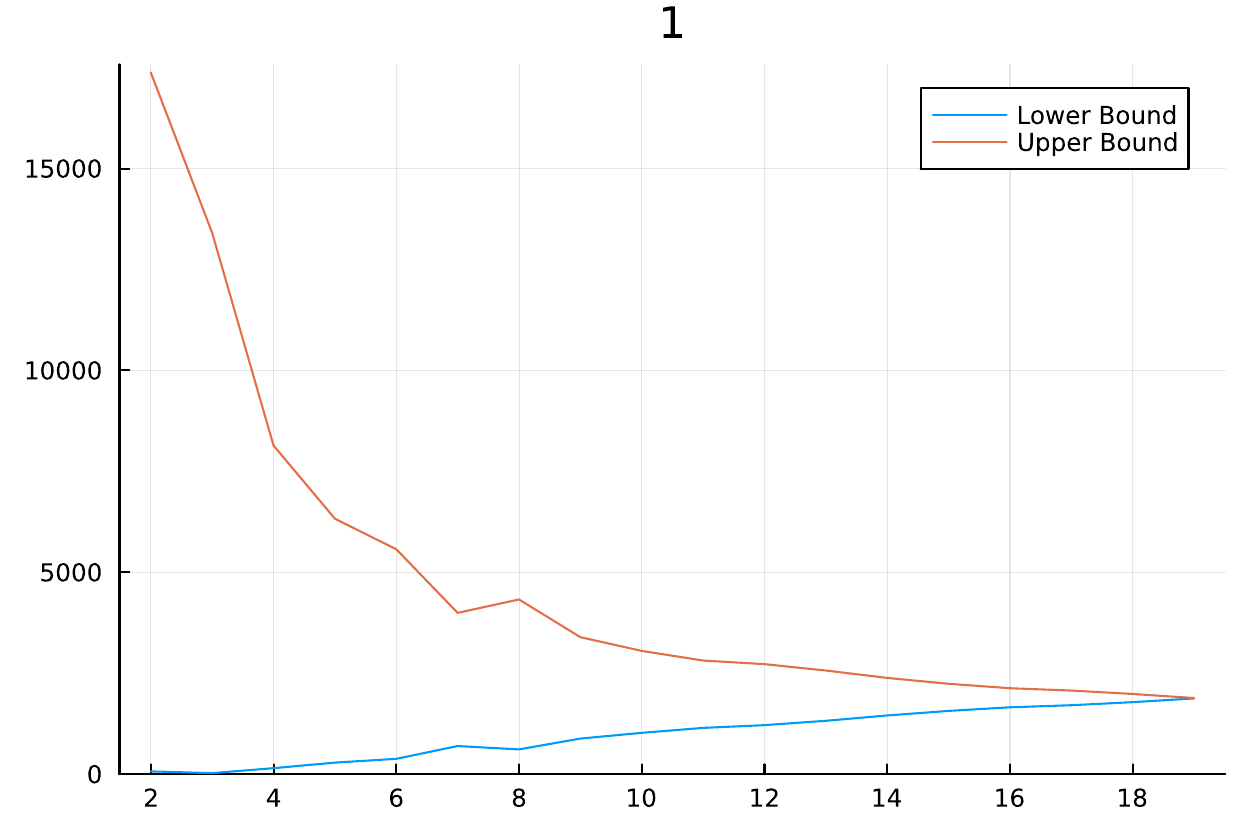}};

  \node[inner sep=0pt, anchor=west] (right) 
    at ([xshift=1.5em]left.east)
    {\includegraphics[width=.42\linewidth]{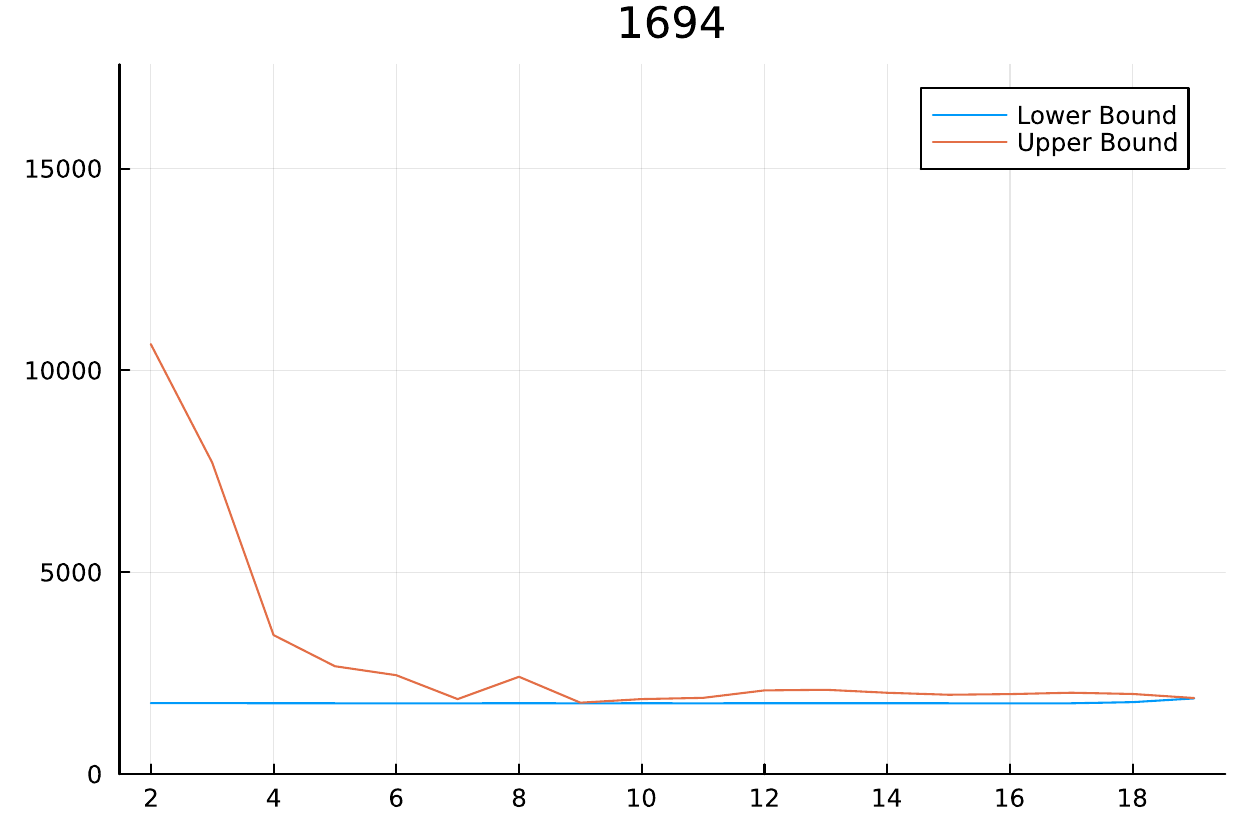}};

  \node[font=\footnotesize, rotate=90, anchor=south]
    at ([xshift=-0.8em]left.west)
    {\color{black}Bounds on LOOCV loss};

  \node[font=\footnotesize, anchor=north]
    at ($(left.south)!0.5!(right.south)+(0,-0.8em)$)
    {\color{black}Cardinality parameter $(\tau)$};
\end{tikzpicture}

\caption{Comparison of initial bounds on LOOCV ($k$-fold with $k=n$) from Algorithm \ref{alg:bounds} (left) and bounds after running Algorithm \ref{alg:parametricK} (right) for a synthetic sparse regression instance where $p=20, n=200, \tau_{\text{true}}=10$, for varying $\tau${\color{black}; see Section \ref{append.syntheticdata} for a full description of our synthetic data generation process}. The black number in the top middle depicts the iteration number of the method.}
\label{fig:algvarykinaction}
\end{figure}


\subsection{Parametric Optimization of $k$-fold Error With Respect to $\gamma$} \label{ssec:parametric2}
In this section, we propose a technique for approximately minimizing the $k$-fold error with respect to the regularization hyperparameter $\gamma$. {\color{black}Note that this corresponds to performing half an iteration of the alternating minimization scheme described at the start of the section, namely approximately minimizing $h$ with respect to $\gamma$ with $\tau$ fixed.}

We begin with two observations from the literature. First, as observed by \cite{stephenson2021can}, the LOOCV error $h(\gamma,\tau)$ is often quasi-convex with respect to $\gamma$ when $\tau=p$. Second, \cite{bertsimas2021unified} {\color{black}and} \cite{bertsimas2022scalable} report that, for sparsity-constrained problems, the optimal support {\color{black}often} does not change as we vary $\gamma$. Combining these observations suggests that, after optimizing $\tau$ with $\gamma$ fixed, a good strategy for minimizing {\color{black}$h$} with respect to $\gamma$ is to fix the optimal support $\bm{z}^{(\mathcal{N}_j)}$ with respect to each fold $\rev{j}$ and invoke a root-finding method to find a $\gamma$ {\color{black}that} locally minimizes 
$\rev{h}$.

Accordingly, we now use the fact that $\gamma$ and $\bm{z}^{(\mathcal{N}_j)}$ fully determine $\bm{\beta}^{(\mathcal{N}_j)}$ to rewrite
\begin{align*}
    \min_{\bm{\beta} \in \mathbb{R}^p} \quad & \frac{\gamma}{2}\Vert \bm{\beta}\Vert_2^2+\Vert \bm{X}\bm{\beta}-\bm{y}\Vert_2^2 \quad \text{s.t.} \quad \beta_i =0 \ \text{if} \ \hat{z}_i=0,\\
\text{as}\qquad &  \bm{\beta}^\star=\left(\frac{\gamma}{2}\mathbb{I}+\bm{X}_{\color{black}\hat{\bm{z}}}^\top \bm{X}_{\color{black}\hat{\bm{z}}}\right)^{-1}\bm{X}_{\color{black}\hat{\bm{z}}}^\top  \bm{y},
\end{align*}
{\color{black}where $\bm{X}_{\hat{\bm{z}}}$ denotes a matrix with the columns of $\bm{X}$ such that $z_i=1$.}

Therefore, we fix each $\bm{z}^{(\mathcal{N}_j)}$ and substitute the resulting expressions for each $\bm{\beta}^{(\mathcal{N}_j)}$ into the $k$-fold error. This substitution yields the following univariate optimization problem, which can be solved via standard root-finding methods to approximately minimize the $k$-fold loss {\color{black}via a local approximation}:
\begin{align}\label{prob:newton1}
            \min_{\gamma >0} \sum_{j \in [k]}\sum_{i \in \mathcal{N}_j} \left(y_i-{\color{black}\bm{x}_i^\top} \mathrm{Diag}(\bm{z}^{(\mathcal{N}_j)}) \left(\frac{\gamma}{2} \mathbb{I}+\mathrm{Diag}(\bm{z}^{(\mathcal{N}_j)})\bm{X}^{(\mathcal{N}_j)\top} \bm{X}^{(\mathcal{N}_j)}\mathrm{Diag}(\bm{z}^{(\mathcal{N}_j)})\right)^{-1}\mathrm{Diag}(\bm{z}^{(\mathcal{N}_j)})\bm{X}^{{(\mathcal{N}_j)}^\top} \bm{y}^{(\mathcal{N}_j)}\right)^2.
\end{align}
Details on minimizing $\gamma$ using \verb|Julia| are provided in \ref{append.4.2}. 

{\color{black}We close this section by remarking that Algorithm \ref{alg:parametricK} could be modified to select the optimal pair $(\tau, \gamma)$ over a grid of candidate values in the same way as it currently selects an optimal $\tau$ with $\gamma$ fixed, although this would likely incur a significantly higher runtime cost than the approach laid out in this paper.}

\section{Numerical Experiments}\label{sec:numres} 

We now present numerical experiments testing our proposed methods. First, in Section~\ref{sec:compExactLOOCV}, we {\color{black}isolate the algorithmic impact of Algorithm~\ref{alg:parametricK} when optimizing the $k$-fold cross-validation error with respect to the sparsity parameter $\tau$. Then, in Section~\ref{sec:comp_real}, we evaluate the statistical performance of the entire \rev{alternating minimization} pipeline proposed in Section~\ref{sec:corddescent} (jointly minimizing $\gamma$ and $\tau$) against widely used sparse regression software packages. {\color{black}We make our code available at \cite{EfficientCrossValidation}.}}

\subsection{Exact {\color{black}$k$}-fold Optimization}\label{sec:compExactLOOCV}

We first assess whether Algorithm~\ref{alg:parametricK} significantly reduces the number of MIOs that need to be solved to minimize the \rev{$k$-fold CV} error with respect to $\tau$, compared to grid search. We set either $k=n$ or $k=10$, corresponding to leave-one-out and 10-fold cross-validation problems \eqref{eq:paramK} respectively.

We compare the performance of two approaches. First, {\color{black}we use} a standard grid search approach (\texttt{Grid}), where we solve the inner MIO in \eqref{eq:paramK} for all combinations of cardinality $\tau\in [p]$ and all folds of the data $j\in [k]$, and select the hyperparameter combination which minimizes the objective. To ensure the quality of the resulting solution, we solve {\color{black}each inner} MIO to optimality (without any time limit{\color{black}, and using the default optimality tolerance of the solver}). Second, we consider using Algorithm~\ref{alg:parametricK} with {\color{black}optimality tolerance $\text{opt}_{\text{tol}}=0.01$ and} parameter $r=\infty$ {\color{black}over the range $\tau \in [p]$} (thus solving MIOs to optimality until the desired optimality gap for problem~\eqref{eq:paramK} is proven). We test regularization parameter $\gamma\in \{0.01,0.02,0.05,0.10,0.20,0.50,1.00\}$ in Algorithm~\ref{alg:parametricK}, and solve all MIOs via their perspective reformulations, namely
\begin{align*}
	\min_{\bm{\beta}\in \mathbb{R}^p,\bm{z}\in \{0,1\}^p}\;& \|\bm{X}\bm{\beta}-\bm{y}\|_2^2+{\color{black}\frac{\gamma-\epsilon}{2}\sum_{j=1}^p \frac{\beta_j^2}{z_j}+\frac{\epsilon}{2}\sum_{i=1}^p\beta_i^2}\quad \text{s.t.}\quad \sum_{j=1}^p z_j\leq \tau,
\end{align*} 
using Mosek 10.0{\color{black}, where $\epsilon=\gamma-0.005$}. Since the approach \texttt{Grid} involves solving $\mathcal{O}(kp)$ MIOs (without a time limit), we are limited to testing these approaches on small datasets, and accordingly use the Diabetes, Housing, Servo, and AutoMPG datasets for this experiment, as described by \cite{gomez2021mixed}. Moreover, we remark that the specific solution times and the number of nodes expanded by each method are not crucial, as those could vary substantially if relaxations other than the perspective are used, different solvers or solution approaches are used, or if advanced techniques are implemented (but both methods would be affected in the same way). Thus, we focus our analysis on relative performance.

Figures~\ref{fig:reductionMioNodes} and~\ref{fig:reductionMioNodes10} summarize the percentage reduction of the number of MIOs and the number of branch-and-bound nodes achieved by Algorithm~\ref{alg:parametricK} over \texttt{Grid}, computed as 
\small$$ \text{Reduction in MIOs}=\frac{\text{\# MIO}_{\texttt{Grid}}-\text{\# MIO}_{\text{Alg.\ref{alg:parametricK}}}}{\text{\# MIO}_{\texttt{Grid}}},\;\; \text{Reduction in nodes}=\frac{\text{\# nodes}_{\texttt{Grid}}-\text{\# nodes}_{\text{Alg.\ref{alg:parametricK}}}}{\text{\# nodes}_{\texttt{Grid}}},$$\normalsize
where $\text{\# MIO}_{\text{Alg.\ref{alg:parametricK}}}$ and $\text{\# nodes}_{\text{Alg.\ref{alg:parametricK}}}$ indicate the number of MIOs or branch-and-bound nodes used by \rev{Algorithm~\ref{alg:parametricK}}.  Tables \ref{tab:gridVsAlg} and \ref{tab:gridVsAlg10} present the detailed computational results. 

\begin{figure}[h!]\centering
\begin{subfigure}[t]{.45\linewidth}
\includegraphics[width=\textwidth,trim={12cm 6cm 12cm 6cm},clip]{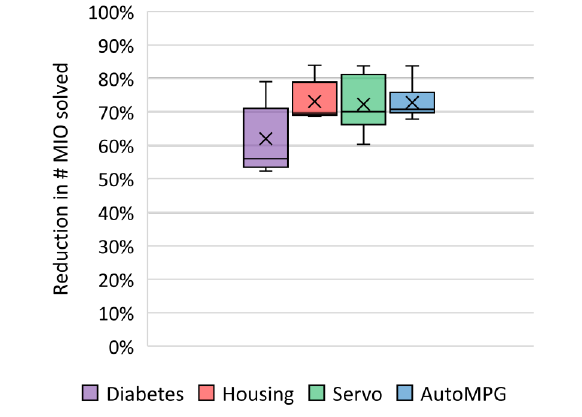}
\end{subfigure}
\begin{subfigure}[t]{.45\linewidth}
\includegraphics[width=\textwidth,trim={12cm 6cm 12cm 6cm},clip]{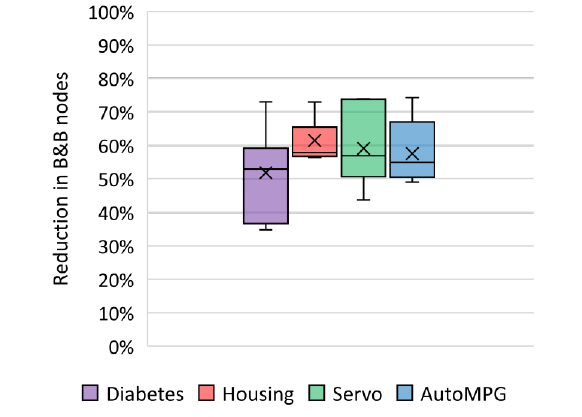}
\end{subfigure}
\caption{Reduction in the number of \rev{MIOs} solved (left) and the total number of branch-and-bound nodes (right) when using Algorithm~\ref{alg:parametricK} for leave-one-out cross-validation, when compared with \texttt{Grid} (i.e., independently solving $\mathcal{O}(pn)$ MIOs) in four real datasets. The distributions shown in the figure correspond to solving the same instance with different values of $\gamma$. All MIOs are solved to optimality, without imposing any time limits. } 
\label{fig:reductionMioNodes}
\end{figure}

\begin{figure}[h!]\centering
\begin{subfigure}[t]{.45\linewidth}
\includegraphics[width=\textwidth,trim={12cm 6cm 12cm 6cm},clip]{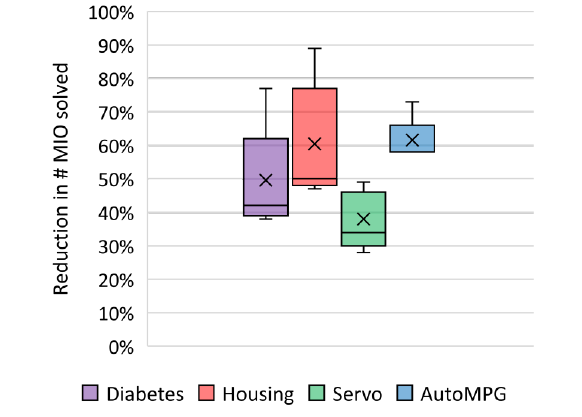}
\end{subfigure}
\begin{subfigure}[t]{.45\linewidth}
\includegraphics[width=\textwidth,trim={12cm 6cm 12cm 6cm},clip]{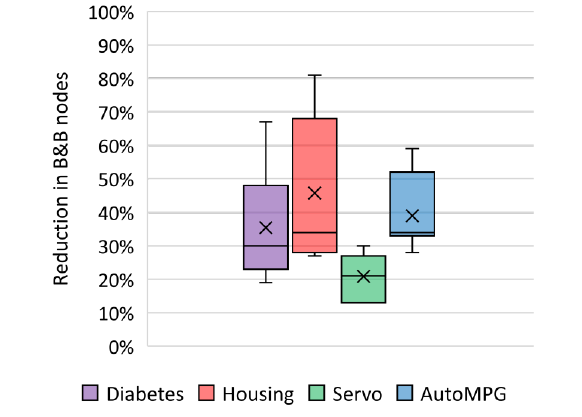}
\end{subfigure}
\caption{\rev{Distribution of the }reduction in the number of MIO{\color{black}s} solved (left) and the total number of branch-and-bound nodes (right) when using Algorithm~\ref{alg:parametricK} for 10-fold cross-validation, when compared with \texttt{Grid} (i.e., independently solving $\mathcal{O}(pk)$ MIOs) in four real datasets \rev{(across different regularization parameters)}. The distributions shown in the figure correspond to solving the same instance with different values of $\gamma$. All MIOs are solved to optimality, without imposing any time limits. } 
\label{fig:reductionMioNodes10}
\end{figure}

	\begin{table}
    [h!]\footnotesize
    \begin{center}
 \caption{Comparison between using Algorithm~\ref{alg:parametricK} and solving $\mathcal{O}(pn)$ MIOs independently (\texttt{Grid}) for leave-one-out cross-validation in four real datasets, for different values of regularization $\gamma$. Times reported are in minutes and correspond to the time to solve all required mixed-integer optimization problems to optimality. No time limits are imposed on the MIOs. Algorithm~\ref{alg:parametricK} consistently reduces {\color{black}the} number of calls to the MIO solver by 50–80\%.}
 \label{tab:gridVsAlg}
	\setlength{\tabcolsep}{2pt}
	\begin{tabular}{l r r r|r r r | r r r  | r r r }
		
		\hline
		\multirow{2}{*}{Dataset}&\multirow{2}{*}{$p$}&\multirow{2}{*}{$n$}&\multirow{2}{*}{$\gamma$}& \multicolumn{3}{c|}{\underline{\texttt{Grid}}}&\multicolumn{3}{c|}{\underline{Algorithm~\ref{alg:parametricK}}}&\multicolumn{3}{c}{\textbf{\underline{Improvement}}}\\
		&&&&Time&\# MIO&Nodes&Time&\# MIO&Nodes&\textbf{Time}&\textbf{\# MIO}&\textbf{Nodes}\\
		\hline
			\multirow{7}{*}{Diabetes}&\multirow{7}{*}{11}&\multirow{7}{*}{442}&
		$0.01$&68&3,978&126,085&37&1,714&59,406&\textbf{45\%}&\textbf{56\%}&\textbf{53\%}\\
		&&&$0.02$&52&3,978&82,523&37&1,768&52,264&\textbf{30\%}&\textbf{56\%}&\textbf{37\%}\\
		&&&$0.05$&42&3,978&42,411&29&1,898&27,652&\textbf{29\%}&\textbf{52\%}&\textbf{35\%}\\
		&&&$0.10$&39&3,978&31,116&26&1,852&16,202&\textbf{34\%}&\textbf{53\%}&\textbf{48\%}\\
		&&&$0.20$&35&3,978&22,165&20&1,332&9,278&\textbf{42\%}&\textbf{67\%}&\textbf{58\%}\\
		&&&$0.50$&32&3,978&11,889&16&1,152&4,852&\textbf{50\%}&\textbf{71\%}&\textbf{59\%}\\
		&&&$1.00$&34&3,978&9,278&14&833&2,501&\textbf{58\%}&\textbf{79\%}&\textbf{73\%}\\
		&&&&&&&&&&&&\\
		\multirow{7}{*}{Housing}&\multirow{7}{*}{13}&\multirow{7}{*}{506}&
		$0.01$&247&6,072&512,723&102&1,906&217,918&\textbf{59\%}&\textbf{69\%}&\textbf{57\%}\\
		&&&$0.02$&187&6,072&324,238&65&1,843&141,493&\textbf{65\%}&\textbf{70\%}&\textbf{56\%}\\
		&&&$0.05$&166&6,072&216,116&92&1,879&93,543&\textbf{45\%}&\textbf{69\%}&\textbf{57\%}\\
		&&&$0.10$&40&6,072&96,387&19&1,880&40,664&\textbf{51\%}&\textbf{69\%}&\textbf{58\%}\\
		&&&$0.20$&82&6,072&68,581&36&1,661&25,171&\textbf{55\%}&\textbf{73\%}&\textbf{63\%}\\
		&&&$0.50$&90&6,072&60,067&34&1,281&20,761&\textbf{62\%}&\textbf{79\%}&\textbf{65\%}\\
		&&&$1.00$&107&6,072&49,770&24&976&13,460&\textbf{77\%}&\textbf{84\%}&\textbf{73\%}\\
		&&&&&&&&&&&&\\
		\multirow{7}{*}{Servo}&\multirow{7}{*}{19}&\multirow{7}{*}{167}&
		$0.01$&466&3,006&1,669,537&276&1,194&940,831&\textbf{41\%}&\textbf{60\%}&\textbf{44\%}\\
		&&&$0.02$&110&3,006&811,432&53&1,016&400,817&\textbf{52\%}&\textbf{66\%}&\textbf{51\%}\\
		&&&$0.05$&44&3,006&324,877&25&986&160,369&\textbf{77\%}&\textbf{84\%}&\textbf{73\%}\\
		&&&$0.10$&23&3,006&162,223&9&686&58,326&\textbf{59\%}&\textbf{77\%}&\textbf{64\%}\\
		&&&$0.20$&15&3,006&76,739&8&900&33,098&\textbf{48\%}&\textbf{70\%}&\textbf{57\%}\\
		&&&$0.50$&10&3,006&40,197&4&566&10,496&\textbf{56\%}&\textbf{81\%}&\textbf{74\%}\\
		&&&$1.00$&8&3,006&25,683&4&488&6,738&\textbf{52\%}&\textbf{84\%}&\textbf{74\%}\\
			&&&&&&&&&&&&\\
		\multirow{7}{*}{AutoMPG}&\multirow{7}{*}{25}&\multirow{7}{*}{392}&
		$0.01$&1,100&9,408&6,772,986&590&3,131&3,532,057&\textbf{46\%}&\textbf{67\%}&\textbf{48\%}\\
		&&&$0.02$&1,356&9,408&3,900,417&450&2,846&1,888,766&\textbf{67\%}&\textbf{70\%}&\textbf{52\%}\\
		&&&$0.05$&519&9,408&2,286,681&227&2,808&1,133,175&\textbf{56\%}&\textbf{70\%}&\textbf{50\%}\\
		&&&$0.10$&355&9,408&1,548,369&145&2,751&687,187&\textbf{59\%}&\textbf{71\%}&\textbf{56\%}\\
		&&&$0.20$&143&9,408&629,020&65&2,686&283,755&\textbf{54\%}&\textbf{71\%}&\textbf{55\%}\\
		&&&$0.50$&66&9,408&176,950&28&2,272&58,464&\textbf{58\%}&\textbf{76\%}&\textbf{67\%}\\
		&&&$1.00$&68&9,408&116,982&38&1,528&30,120&\textbf{43\%}&\textbf{84\%}&\textbf{74\%}\\
		\hline
	\end{tabular}
    \end{center}
\end{table}

	\begin{table}[h!]\footnotesize
    \begin{center}
 \caption{Comparison between using Algorithm~\ref{alg:parametricK} and solving $\mathcal{O}(pk)$ MIOs independently (\texttt{Grid}) for 10-fold cross validation in four real datasets, for different values of regularization $\gamma$. Times reported are in minutes, and correspond to the time to solve all required mixed-integer optimization problems to optimality. No time limits are imposed on the MIOs.}
 \label{tab:gridVsAlg10}
	\setlength{\tabcolsep}{2pt}
	\begin{tabular}{l r r r|r r r | r r r  | r r r }
		
		\hline
		\multirow{2}{*}{Dataset}&\multirow{2}{*}{$p$}&\multirow{2}{*}{$n$}&\multirow{2}{*}{$\gamma$}& \multicolumn{3}{c|}{\underline{\texttt{Grid}}}&\multicolumn{3}{c|}{\underline{Algorithm~\ref{alg:parametricK}}}&\multicolumn{3}{c}{\textbf{\underline{Improvement}}}\\
		&&&&Time&\# MIO&Nodes&Time&\# MIO&Nodes&\textbf{Time}&\textbf{\# MIO}&\textbf{Nodes}\\
		\hline
			\multirow{7}{*}{Diabetes}&\multirow{7}{*}{11}&\multirow{7}{*}{442}&
		$0.01$&3&396&11,666&2&242&8,224&\textbf{14\%}&\textbf{39\%}&\textbf{30\%}\\
		&&&$0.02$&2&396&8,371&2&235&6,785&\textbf{12\%}&\textbf{41\%}&\textbf{19\%}\\
		&&&$0.05$&2&396&4,436&2&228&3,430&\textbf{10\%}&\textbf{42\%}&\textbf{23\%}\\
		&&&$0.10$&2&396&3,185&2&247&2,277&\textbf{10\%}&\textbf{38\%}&\textbf{29\%}\\
		&&&$0.20$&1&396&2,268&1&206&1,536&\textbf{8\%}&\textbf{48\%}&\textbf{32\%}\\
		&&&$0.50$&1&396&1,233&1&149&643&\textbf{26\%}&\textbf{62\%}&\textbf{48\%}\\
		&&&$1.00$&1&396&872&1&93&287&\textbf{42\%}&\textbf{77\%}&\textbf{67\%}\\
		&&&&&&&&&&&&\\
		\multirow{7}{*}{Housing}&\multirow{7}{*}{13}&\multirow{7}{*}{506}&
		$0.01$&25&600&48,069&19&321&35,227&\textbf{25\%}&\textbf{47\%}&\textbf{27\%}\\
		&&&$0.02$&19&600&34,915&14&310&25,090&\textbf{28\%}&\textbf{48\%}&\textbf{28\%}\\
		&&&$0.05$&14&600&21,350&10&303&14,933&\textbf{29\%}&\textbf{50\%}&\textbf{30\%}\\
		&&&$0.10$&10&600&11,012&7&300&7,308&\textbf{31\%}&\textbf{50\%}&\textbf{34\%}\\
		&&&$0.20$&9&600&7,406&5&230&3,524&\textbf{46\%}&\textbf{62\%}&\textbf{52\%}\\
		&&&$0.50$&9&600&6,168&3&141&1,977&\textbf{62\%}&\textbf{77\%}&\textbf{68\%}\\
		&&&$1.00$&8&600&4,993&2&66&930&\textbf{77\%}&\textbf{89\%}&\textbf{81\%}\\
		&&&&&&&&&&&&\\
		\multirow{7}{*}{Servo}&\multirow{7}{*}{19}&\multirow{7}{*}{167}&
		$0.01$&15&288&148,168&12&191&128,592&\textbf{16\%}&\textbf{34\%}&\textbf{13\%}\\
		&&&$0.02$&8&288&77,457&7&190&67,416&\textbf{10\%}&\textbf{34\%}&\textbf{13\%}\\
		&&&$0.05$&3&288&29,056&3&157&23,653&\textbf{16\%}&\textbf{45\%}&\textbf{19\%}\\
		&&&$0.10$&2&288&15,951&2&146&12,562&\textbf{16\%}&\textbf{49\%}&\textbf{21\%}\\
		&&&$0.20$&1&288&8,117&1&155&6,275&\textbf{12\%}&\textbf{46\%}&\textbf{23\%}\\
		&&&$0.50$&1&288&4,028&1&201&2,922&\textbf{3\%}&\textbf{30\%}&\textbf{27\%}\\
		&&&$1.00$&1&288&2,541&1&206&1,768&\textbf{1\%}&\textbf{28\%}&\textbf{30\%}\\
			&&&&&&&&&&&&\\
		\multirow{7}{*}{AutoMPG}&\multirow{7}{*}{25}&\multirow{7}{*}{392}&
		$0.01$&111&936&691,816&76&389&460,187&\textbf{31\%}&\textbf{58\%}&\textbf{33\%}\\
		&&&$0.02$&68&936&401,905&44&374&264,179&\textbf{35\%}&\textbf{60\%}&\textbf{34\%}\\
		&&&$0.05$&42&936&225,318&30&396&161,639&\textbf{28\%}&\textbf{58\%}&\textbf{28\%}\\
		&&&$0.10$&30&936&149,243&20&389&98,261&\textbf{35\%}&\textbf{58\%}&\textbf{34\%}\\
		&&&$0.20$&14&936&61,534&10&389&41,323&\textbf{32\%}&\textbf{58\%}&\textbf{33\%}\\
		&&&$0.50$&7&936&17,865&4&318&8,550&\textbf{43\%}&\textbf{66\%}&\textbf{52\%}\\
		&&&$1.00$&6&936&10,848&3&251&4,480&\textbf{48\%}&\textbf{73\%}&\textbf{59\%}\\
		\hline
	\end{tabular}
    \end{center}
\end{table}

We observe that across these four datasets, Algorithm~\ref{alg:parametricK} reduces the number of \rev{MIOs} that need to be solved by an average of 70\% for leave-one-out cross-validation and by 52\% for 10-fold cross-validation. The overall number of branch-and-bound nodes is reduced by an average of 57\% for leave-one-out cross-validation and 35\% for 10-fold cross-validation (the reduction in computational times is similar to the reduction of nodes). \rev{Note that MIOs with strong continuous relaxations are less likely to be solved to optimality by Algorithm~\ref{alg:parametricK}. In general, these MIOs are also easier to solve to optimality, thus resulting in smaller improvements in {\color{black}the} number of nodes and solution times than those suggested by the number of MIOs solved. Nonetheless, this discrepancy is relatively small, indicating that the proposed approach is still effective in avoid{\color{black}ing} solving several non-trivial MIOs.}

We observe that solution times for both methods decrease on a given dataset as $\gamma$ increases (as expected, since the perspective reformulation is stronger). Interestingly, while the improvements of Algorithm~\ref{alg:parametricK} over \texttt{Grid} (in terms of time, MIOs solved, and nodes) are more pronounced in regimes with large regularization $\gamma$, this effect on $\gamma$ is slight: Algorithm~\ref{alg:parametricK} consistently results in improvements over 40\% (and often more) even for the smallest values of $\gamma$ tested. These results indicate that the relaxations of the bilevel optimization \eqref{eq:paramK} derived in \S\ref{sec:lowerbounds} are sufficiently strong to avoid solving most of the MIOs that traditional methods such as \texttt{Grid} would solve, without sacrificing solution quality. The proposed methods are especially beneficial for settings where $k$ is large, that is, in the settings that would require more MIOs and are more computationally expensive using standard approaches. 

The resulting approach still requires solving several MIOs, but, as we show throughout the rest of this section, approximating each MIO with its perspective relaxation yields similarly high-quality statistical estimators at a fraction of the computational cost.

\subsection{Statistical Results With Real Data}\label{sec:comp_real}

{\color{black}In this section, we investigate the performance of our alternating minimization procedure as described in Section \ref{sec:corddescent}  (called MIO throughout the section) on a suite of eleven real datasets from the UCI repository. Specifically, we use the following experimental setup for the MIO method:
\begin{itemize}
    \item We use a grid of ten different values of $\gamma$ log-uniformly distributed on $[10^{-4}, 10^4]$, with $\gamma_0=\frac{1}{\sqrt{n}}$.
    \item We set $\tau_{\min}=2$ {\color{black}and} $\tau_{\max}$ {\color{black}to be the largest integer }such that $\tau_{\max}\log \tau_{\max}\leq \min(n,p)$ when optimizing $\tau$ in this experiment\footnote{We previously tried setting $\tau_{\max}=p$. We found that this yielded the same optimal hyperparameters, but increased the total runtime of the method substantially.} as in \cite{gamarnik2022sparse}.
    \item We impose a limit of at most $10$ iterations of alternating minimization on minimizing the k-fold cross-validation error with respect to $\gamma$ and with respect to $\tau$ (terminating early if $(\gamma_t, \tau_t)=(\gamma_{t-1}, \tau_{t-1})$).
    \item For tractability, for the overdetermined (underdetermined) instances, we impose a time limit of $120$s ($600$s) for each MIO solved when cross-validating $\tau$, and a time limit of $3600$s ($7200$s) for fitting the final MIO once all hyperparameters are fixed. When we exceed a time limit, we use the best regression model found by the MIO solver at that time.
\end{itemize}
 We refer to this implementation of our \rev{alternating minimization} approach as ``MIO'' (short for mixed-integer optimization).}

We compare against the following state-of-the-art methods, using \rev{built-in} functions to approximately minimize the cross-validation loss with respect to the method's hyperparameters via grid search, and subsequently fit a regression model on the entire dataset with these cross-validated parameters (see also \cite{bertsimas2020sparse2} for a detailed discussion of these approaches):
\begin{itemize}
    \item The \verb|ElasticNet| method in the ubiquitous \verb|glmnet| package, with grid search on their parameter $\alpha \in \{0, 0.1, 0.2, \ldots, 1\}$
    \item The Minimax Concave Penalty (MCP) 
    as implemented in the \verb|R| package \verb|ncvreg|, using the \verb|cv.ncvreg| function with 
    default parameters to minimize the {\color{black}five-fold} cross-validation error.\color{black}
    \item The \verb|L0Learn.cvfit| method implemented in the \verb|L0Learn| \verb|R| package \citep[cf.][]{hazimeh2020fast}, with five folds, a grid of $10$ different values of $\gamma$ and default parameters otherwise.
\end{itemize}
{\color{black} Note that we use default parameters for glmnet, MCP, and L0Learn to reflect typical practitioner usage. In particular, we do not impose an explicit cardinality budget for ElasticNet or MCP because neither method provides a parameter option to do so. Similarly, we do not explicitly impose a cardinality constraint in \verb|L0Learn|, because \verb|L0Learn| penalizes rather than constrains cardinality.}

We compare performance in terms of the Mean Square{\color{black}d} Error, namely
\begin{align*}
    MSE(\bm{\beta}):=\frac{1}{n}\sum_{i=1}^n (y_{i}-\bm{x}_{i}^\top \bm{\beta})^2,
\end{align*}
which can either be taken over the validation set (CV)--that is, the objective \eqref{prob:upperlevelofv_validation} we attempt to minimize--or over an unseen test set (MSE), acting as a proxy for generalization error. 

To measure the validation and test set errors, we repeat the following procedure {\color{black}five} times and report the average result: we randomly shuffle the data into $80\%$ training/validation and $20\%$ test data, perform five-fold cross-validation on the $80\%$ training/validation data, fit a model with the cross-validated $(\tau, \gamma)$ on the combined $80\%$ train/validation data, and evaluate the model's test-set performance on the remaining $20\%$ test data. We also report the average value of $\tau$, the cross-validated sparsity, for each method. {\color{black}Note that the use of a ``2'' after the dataset name indicates that a dataset includes second-order interactions, in order to increase the computational difficulty of processing the dataset.}


{\color{black}
We observe in Tables \ref{tab:comparison_ourmethods5}--\ref{tab:comparison_ourmethods6} that MIO performs comparably to widely used methods for most datasets, especially when they are relatively underdetermined. In particular, across the five most underdetermined datasets considered, it yields a five-fold cross-validation error $21.2\%$ lower than MCP, $7.3\%$ lower than glmnet, and $1.5\%$ lower than L0Learn. However, it admittedly yields a higher average cross-validation error on the more overdetermined datasets (by $14.6\%$, $37.6\%$ and $27.4\%$ respectively). In terms of test-set errors, for the five most underdetermined datasets, MIO performs $6.9\%$ better than MCP, $2.1\%$ better than glmnet, and $4.2\%$ better than L0Learn on average. However, it performs $1.9\%$, $18.2\%$, and $19.1\%$ worse on the more overdetermined datasets in terms of test-set error, respectively. 
All in all, our results show that MIO-based cross-validation performs best on relatively underdetermined problems, whereas for highly overdetermined datasets, standard packages such as glmnet or L0Learn remain preferable from both statistical and computational perspectives.

Finally, we note that the average runtime across all datasets was $47510$s for MIO (median: $69.90$s, range: $0.42$s--$541511$s), $0.180$s for MCP, $1.154$s for glmnet, and $0.841$s for L0Learn, respectively. This shows that while our scheme is significantly slower than state-of-the-art regression packages with sophisticated codebases, our cross-validation scheme does return reasonable results on real-world datasets. 
}

\begin{table}[h!]
\centering\footnotesize
\caption{\color{black} Average ($\pm$ one standard deviation) performance of five-fold versions of the methods across a suite of real-world datasets where the ground truth is unknown (and may not be sparse), sorted by how overdetermined the dataset is ($n/p$), and separated into the underdetermined and overdetermined cases, where $\tau:=\Vert \bm{\beta}\Vert_0$ denotes the sparsity of the regression model $\bm{\beta}$. In underdetermined settings, MIO often yields competitive or lower out-of-sample MSEs than MCP, glmnet, and L0Learn, whereas in more overdetermined settings, glmnet or L0Learn tend to yield lower out-of-sample MSEs.
}
\begin{tabular}{@{}l r r r r r r r r @{}} \toprule
Dataset & n & p & \multicolumn{3}{c@{\hspace{0mm}}}{MIO} &  \multicolumn{3}{c@{\hspace{0mm}}}{MCP} \\ 
\cmidrule(l){4-6} \cmidrule(l){7-9}  &   &   & $\tau$ & CV & MSE&$\tau$ & CV & MSE \\\midrule
Wine & 6497 & 11 & 6 $\pm$ 0 & 0.567 $\pm$ 0.003 & 0.545 $\pm$ 0.021 & 10.8 $\pm$ 0.447 & 0.543 $\pm$ 0.005 & 0.543 $\pm$ 0.020 \\
AutoMPG & 392 & 25 & 11 $\pm$ 0 & 10.192 $\pm$ 0.722 & 9.595 $\pm$ 2.901 & 16.6 $\pm$ 1.817 & 9.079 $\pm$ 0.486 & 9.043 $\pm$ 1.970 \\
Hitters & 263 & 19 & 8.6 $\pm$ 0.894 & 0.077 $\pm$ 0.006 & 0.079 $\pm$ 0.02 & 12.8 $\pm$ 4.604 & 0.08 $\pm$ 0.008 & 0.081 $\pm$ 0.024 \\
Prostate & 97 & 8 & 3.8 $\pm$ 1.304 & 0.53 $\pm$ 0.035 & 0.609 $\pm$ 0.159 & 5.8 $\pm$ 2.168 & 0.572 $\pm$ 0.051 & 0.574 $\pm$ 0.153 \\
Servo & 167 & 19 & 7.6 $\pm$ 2.608 & 0.785 $\pm$ 0.106 & 0.769 $\pm$ 0.173 & 14.000 $\pm$ 1.000 & 0.752 $\pm$ 0.062 & 0.723 $\pm$ 0.212 \\
Housing2 & 506 & 91 & 22.6 $\pm$ 3.05 & 19.182 $\pm$ 1.57 & 16.776 $\pm$ 7.548 & 34.2 $\pm$ 5.718 & 16.311 $\pm$ 2.325 & 16.879 $\pm$ 3.449 \\
Toxicity & 38 & 9 & 3.6 $\pm$ 1.342 & 0.036 $\pm$ 0.009 & 0.05 $\pm$ 0.037 & 2.6 $\pm$ 0.894 & 0.049 $\pm$ 0.01 & 0.057 $\pm$ 0.049 \\
Steam & 25 & 8 & 3.4 $\pm$ 1.14 & 0.463 $\pm$ 0.084 & 0.490 $\pm$ 0.300 & 2.6 $\pm$ 0.894 & 0.629 $\pm$ 0.103 & 0.532 $\pm$ 0.200 \\
Alcohol2 & 44 & 21 & 3.6 $\pm$ 1.517 & 0.225 $\pm$ 0.032 & 0.256 $\pm$ 0.094 & 1.8 $\pm$ 0.447 & 0.238 $\pm$ 0.035 & 0.258 $\pm$ 0.047 \\ \midrule
TopGear & 242 & 373 & 18.2 $\pm$ 10.134 & 0.046 $\pm$ 0.006 & 0.049 $\pm$ 0.005 & 7.4 $\pm$ 1.517 & 0.061 $\pm$ 0.006 & 0.061 $\pm$ 0.018 \\ 
Bardet & 120 & 200 & 18.8 $\pm$ 5.891 & 0.007 $\pm$ 0 & 0.009 $\pm$ 0.003 & 5.6 $\pm$ 1.14 & 0.009 $\pm$ 0.002 & 0.009 $\pm$ 0.002 \\
\bottomrule
\end{tabular}

\label{tab:comparison_ourmethods5}
\end{table}

\begin{table}[h!]
\centering\footnotesize
\caption{\color{black} Average ($\pm$ one standard deviation) performance of five-fold versions of the methods across a suite of real-world datasets where the ground truth is unknown (and may not be sparse), sorted by how overdetermined the dataset is ($n/p$), and separated into the underdetermined and overdetermined cases (cont).
}
\begin{tabular}{@{}l r r r r r r r r r@{}} \toprule
Dataset & n & p & \multicolumn{3}{c@{\hspace{0mm}}}{glmnet}& \multicolumn{3}{c@{\hspace{0mm}}}{L0Learn} \\ 
\cmidrule(l){4-6} \cmidrule(l){7-9} &   &   & $\tau$ & CV & MSE&$\tau$ & CV & MSE \\ \midrule
Wine & 6497 & 11 & 11 $\pm$ 0 & 0.542 $\pm$ 0.005 & 0.543 $\pm$ 0.02 & 10.6 $\pm$ 0.548 & 0.542 $\pm$ 0.005 & 0.543 $\pm$ 0.02\\
AutoMPG & 392 & 25 & 22.6 $\pm$ 1.949 & 8.627 $\pm$ 0.493 & 9.201 $\pm$ 2.436 & 15.8 $\pm$ 4.266 & 9.099 $\pm$ 0.616 & 9.061 $\pm$ 2.351\\
Hitters & 263 & 19 & 14.4 $\pm$ 4.561 & 0.077 $\pm$ 0.006 & 0.082 $\pm$ 0.023 & 10.6 $\pm$ 5.55 & 0.075 $\pm$ 0.006 & 0.08 $\pm$ 0.021\\
Prostate & 97 & 8 & 6.8 $\pm$ 1.095 & 0.507 $\pm$ 0.051 & 0.581 $\pm$ 0.154 & 5.6 $\pm$ 2.408 & 0.508 $\pm$ 0.048 & 0.569 $\pm$ 0.183\\
Servo & 167 & 19 & 16 $\pm$ 1.414 & 0.693 $\pm$ 0.051 & 0.726 $\pm$ 0.211 & 10.2 $\pm$ 3.493 & 0.696 $\pm$ 0.079 & 0.746 $\pm$ 0.221 \\
Housing2 & 506 & 91 & 86 $\pm$ 2.915 & 12.317 $\pm$ 0.282 & 12.867 $\pm$ 2.133 & 58.6 $\pm$ 9.099 & 13.669 $\pm$ 0.997 & 12.818 $\pm$ 2.009\\
Toxicity & 38 & 9 & 6 $\pm$ 0.707 & 0.041 $\pm$ 0.009 & 0.047 $\pm$ 0.029 & 3 $\pm$ 1.871 & 0.037 $\pm$ 0.012 & 0.062 $\pm$ 0.049\\
Steam & 25 & 8 & 4.4 $\pm$ 0.894 & 0.492 $\pm$ 0.14 & 0.507 $\pm$ 0.200 & 2.2 $\pm$ 0.447 & 0.475 $\pm$ 0.179 & 0.499 $\pm$ 0.103\\
Alcohol2 & 44 & 21 & 8.8 $\pm$ 4.087 & 0.253 $\pm$ 0.05 & 0.263 $\pm$ 0.058 & 7.6 $\pm$ 6.269 & 0.220 $\pm$ 0.021 & 0.272 $\pm$ 0.072 \\\midrule
TopGear & 242 & 373 & 39.4 $\pm$ 21.686 & 0.045 $\pm$ 0.004 & 0.047 $\pm$ 0.014 & 28.8 $\pm$ 37.626 & 0.050 $\pm$ 0.010 & 0.049 $\pm$ 0.015\\
Bardet & 120 & 200 & 29.8 $\pm$ 6.907 & 0.007 $\pm$ 0.001 & 0.008 $\pm$ 0.003 & 30.4 $\pm$ 13.446 & 0.007 $\pm$ 0.000 & 0.009 $\pm$ 0.004\\
\bottomrule
\end{tabular}

\label{tab:comparison_ourmethods6}
\end{table}

\section{Conclusion}
In this paper, we propose a new optimization-based approach 
for {\color{black}selecting hyperparameters via cross-validation} in ridge-regularized sparse regression problems, by leveraging perspective relaxations and bounds on the cross-validation error. The proposed approach substantially decreases the number of MIOs and branch-and-bound nodes explored to optimize the cross-validation loss. {\color{black}Overall, these results suggest that perspective relaxations can help to make MIO‑driven cross‑validation practically viable. 
As future work, it could be interesting to explore how strong convex relaxations can help accelerate MIO-driven cross-validation in other contexts, e.g., when designing optimal decision trees or neural networks.}

 \subsubsection*{Acknowledgments: } Andr\'es G\'omez is supported in part by grant FA9550-24-1-0086 from the Air Force Office of Scientific Research. Ryan Cory-Wright gratefully acknowledges the MIT-IBM Research Lab for hosting him while part of this work was conducted. We are grateful to Jean Pauphilet for valuable discussions on an earlier draft of this manuscript.
{\def\enotesize{\footnotesize}

\theendnotes
}

\FloatBarrier
{
\scriptsize
\bibliographystyle{informs2014}

}

\newpage
\ECSwitch
\ECHead{Supplementary Material}

\section{Implementation Details}\label{sec:append.implementationdetails}
To solve each MIO in Algorithm \ref{alg:parametricK}, we invoke a Generalized Benders Decomposition scheme \citep{geoffrion1972generalized}, which was specialized to sparse regression problems by \cite{bertsimas2020sparse}. For any fixed $\gamma, \tau$, the method proceeds by minimizing a piecewise linear approximation of 
\begin{align}
    f(\bm{z}, \gamma):=\min_{\bm{\beta} \in \mathbb{R}^p} \ \frac{\gamma}{2}\sum_{j \in [p]}\frac{\beta_j^2}{z_j}+\Vert \bm{X}^{}\bm{\beta}-\bm{y}^{}\Vert_2^2,\label{prob:persp1}
\end{align}
{\color{black}with respect to the support $\bm{z}$,}
until it converges to an optimal solution or encounters a time limit. 

We now discuss two enhancements that improve this method's performance in practice.

\paragraph{Warm-Starts:} as noted by \cite{ bertsimas2021unified}, a greedily rounded solution to the Boolean relaxation constitutes an excellent warm-start for a Generalized Benders Decomposition scheme. Therefore, when computing the lower and upper bounds on $h_{j}(\gamma, \tau)$ for each $\tau$ by solving a perspective relaxation, we save the greedily rounded solution to the relaxation in memory, and provide the relevant rounding as a high-quality warm-start before solving the corresponding MIO. 

\paragraph{Screening Rules:} as observed by \cite{atamturk2020safe}, if we have an upper bound on the optimal value of $f(\bm{z}, \gamma)$, say $\bar{f}$, 
an optimal solution to the Boolean relaxation of minimizing \eqref{prob:persp1} over $\bm{z} \in [0, 1]^p$, say $(\bm{\beta}, \bm{z})$, and a lower bound on the optimal value of $h(\bm{z}, \gamma)$ from the Boolean relaxation, say $\underaccent{\bar}{f}$, then, letting $\beta_{[\tau]}$ be the $\tau$th largest value of $\beta$ in absolute magnitude, we have the following screening rules:
    \begin{itemize}
        \item If $\beta_i^2 \leq \beta_{[\tau+1]}^2$ and $\underaccent{\bar}{f}-\frac{1}{2\gamma}(\beta^2_{i}-\beta^2_{[\tau]})>\bar{f}$ then $z_i=0$.
        \item If $\beta_i^2 \geq \beta_{[\tau]}^2$ and $\underaccent{\bar}{f}+\frac{1}{2\gamma}(\beta_{i}^2-\beta^2_{[\tau+1]})>\bar{f}$ then $z_i=1$.
    \end{itemize}
Accordingly, to reduce the dimensionality of our problems, we solve a perspective relaxation for each fold of the data with $\tau=\tau_{\max}$ as a preprocessing step, and screen out the features where $z_i=0$ at $\tau=\tau_{\max}$ (for this fold of the data) before running Generalized Benders Decomposition. 

{\color{black}As reported by \citet{atamturk2020safe}, screening rules often reduce the number of decision variables in an MIO by $20\%$--$97\%$, with the most significant benefits when $\gamma$ is relatively large, and the duality gap between an MIO and its perspective relaxation is relatively small.}

\subsection{Implementation Details for Section \ref{ssec:parametric2}}\label{append.4.2}
{\color{black}
In our numerical experiments, we find local minimizers of our approximation of {\color{black}$h$} by invoking the \verb|ForwardDiff| function in \verb|Julia| to automatically differentiate our approximation of {\color{black}$h$}, and subsequently identify local minima via the \verb|Order0| method in the \verb|Roots.jl| package, which is designed to be a robust root-finding method. To avoid convergence to a low-quality local minimum, we run the search algorithm initialized at the previous iterate $\gamma_{t-1}$ and ten points log-uniformly distributed in $[10^{-4}, 10^{4}]$, and set $\gamma_t$ to be a local minimum with the smallest estimated error. Moreover, to ensure numerical robustness, we require that $\gamma_t$ remains within the bounds $[10^{-4}, 10^4]$ and project $\gamma_t$ onto this interval if it exceeds these bounds (this almost never occurs in practice, because the data is preprocessed to be standardized). This approach is highly efficient in practice, particularly when the optimal support does not vary significantly with $\gamma$. 
}

\subsection{Synthetic Data Generation Process}\label{append.syntheticdata}
{\color{black}For our synthetic experiments, we follow the experimental setup in \cite{bertsimas2020sparse2}. Given a fixed number of features $p$, number of data points $n$, true sparsity $1\leq \tau_{\text{true}}\leq p$, autocorrelation parameter $0\leq \rho\leq 1$ and signal-to-noise parameter $\nu$:
\begin{enumerate}
\item The rows of the model matrix are generated i.i.d. from a $p$-dimensional multivariate Gaussian distribution $\mathcal{N}(\bm{0},\bm{\Sigma})$, where $\Sigma_{ij}=\rho^{|i-j|}$ for all $i,j\in [p]$.
\item A ``ground-truth" vector $\bm{\beta}_{\text{true}}$ is sampled with exactly $\tau_{\text{true}}$ non-zero coefficients. The position of the non-zero entries is randomly chosen from a uniform distribution, and the value of the non-zero entries is either $1$ or $-1$ with equal probability.
\item The response vector is generated as $\bm{y}=\bm{X\beta_{\text{true}}}+\bm{\varepsilon}$, where each $\varepsilon_i$ is generated iid from a scaled normal distribution such that $\sqrt{\nu}=\|\bm{X\beta_{\text{true}}}\|_2/\|\bm{\varepsilon}\|_2$. 
\item We standardize $\bm{X}, \bm{y}$ to normalize and center them.
\end{enumerate}
}

{\color{black}
\section{Non-dominance of Relaxations  \eqref{eqn:perspbounds}  and \eqref{eq:LOOLower}}\label{ec.nondominanceofrelaxations}
In this section, we construct examples demonstrating that the lower bounds in \eqref{eqn:perspbounds} and \eqref{eq:LOOLower} are non-dominated by each other.

\textbf{Problem \eqref{eqn:perspbounds} can be stronger than Problem \eqref{eq:LOOLower}: } Consider a setting with
$$
n = 3,\ p = 1,\ X = \begin{pmatrix}1 \\ 2\\ 3\end{pmatrix},\quad
y = \begin{pmatrix}1\\ 2\\ 5\\\end{pmatrix}, k=3
$$
and hyperparameters
$
\gamma = 1,\ \tau = 1.
$ 

Then, since the sparsity constraint is not binding (for either the problem of training over the entire dataset or with one fold left out), our lower-level problems all reduce to ridge regression in one dimension, and our perspective relaxations are all tight. In particular, the full-data training problem admits an optimal solution $\beta^\star=40/29$ with an optimal objective value of $2.414$. On the other hand, on leaving out the third observation and training on the first two folds, we obtain $\beta^\star=10/11$ with training objective $0.455$, with a fold-$3$ contribution to the LOOCV error of $5.165$. In particular, the lower bound from \eqref{eqn:perspbounds} is $625/121$, which is larger than the lower bound from \eqref{eq:LOOLower}, namely $625/319$. 

\textbf{Problem \eqref{eqn:perspbounds} can be weaker than Problem \eqref{eq:LOOLower}: }  
Consider a setting with
$$
n = 3,\ p = 2,\ X = \begin{pmatrix}1 & 1\\ -2 & 0\\ 2 &1\end{pmatrix},\quad
y = \begin{pmatrix}0\\ 3\\ 3\\\end{pmatrix}, k=3
$$
and hyperparameters
$
\gamma = 1,\ \tau = 1,
$  where we focus on leaving out the first fold. Then, the optimal full-data solution is $\beta^\star=(0, 6/5)$ with objective value $72/5$. On the other hand, if we leave out the first observation then the optimal solution becomes $\beta^\star=(0, 2)$ with objective value $12$, and optimal value of the perspective relaxation of $10.714$. Accordingly, the lower bound on the LOOCV error of fold $1$ from \eqref{eq:LOOLower} is $2.4$, while the lower bound from \eqref{eqn:perspbounds} is $9/7 < 2.4$. Thus, \eqref{eqn:perspbounds} can be weaker than \eqref{eq:LOOLower} even when using the optimal value of the MIO and perspective relaxation for upper and lower bounds in \eqref{eqn:perspbounds}.

}
\end{document}